%% file: PopulationBalance.tex
\documentclass[final,preprint,3p,sort&compress]{elsarticle}
\journal{{\tt arXiv.org}}
\usepackage[dvips]{epsfig}
\usepackage{graphicx} 
\usepackage{graphics}
\usepackage{pdfpages}
\usepackage{floatrow}
\usepackage{latexsym}
\usepackage{verbatim}
\usepackage{amsmath,amsthm}
\usepackage{subfig}
\usepackage{amsthm}
\usepackage{amssymb}
\usepackage{float}
\usepackage{caption}
\usepackage{bbm}
\usepackage{diagbox}

\usepackage{fonts}
\usepackage{enumerate}
\usepackage{mathrsfs}
\usepackage{placeins}
\usepackage{standalone}
\usepackage{paralist}

\definecolor{greenyellow}   {cmyk}{0.15, 0   , 0.69, 0   }
\definecolor{yellow}        {cmyk}{0   , 0   , 1   , 0   }
\definecolor{goldenrod}     {cmyk}{0   , 0.10, 0.84, 0   }
\definecolor{dandelion}     {cmyk}{0   , 0.29, 0.84, 0   }
\definecolor{apricot}       {cmyk}{0   , 0.32, 0.52, 0   }
\definecolor{peach}         {cmyk}{0   , 0.50, 0.70, 0   }
\definecolor{melon}         {cmyk}{0   , 0.46, 0.50, 0   }
\definecolor{yelloworange}  {cmyk}{0   , 0.42, 1   , 0   }
\definecolor{orange}        {cmyk}{0   , 0.61, 0.87, 0   }
\definecolor{burntorange}   {cmyk}{0   , 0.51, 1   , 0   }
\definecolor{bittersweet}   {cmyk}{0   , 0.75, 1   , 0.24}
\definecolor{redorange}     {cmyk}{0   , 0.77, 0.87, 0   }
\definecolor{mahogany}      {cmyk}{0   , 0.85, 0.87, 0.35}
\definecolor{maroon}        {cmyk}{0   , 0.87, 0.68, 0.32}
\definecolor{brickred}      {cmyk}{0   , 0.89, 0.94, 0.28}
\definecolor{red}           {cmyk}{0   , 1   , 1   , 0   }
\definecolor{orangered}     {cmyk}{0   , 1   , 0.50, 0   }
\definecolor{rubinered}     {cmyk}{0   , 1   , 0.13, 0   }
\definecolor{wildstrawberry}{cmyk}{0   , 0.96, 0.39, 0   }
\definecolor{salmon}        {cmyk}{0   , 0.53, 0.38, 0   }
\definecolor{carnationpink} {cmyk}{0   , 0.63, 0   , 0   }
\definecolor{magenta}       {cmyk}{0   , 1   , 0   , 0   }
\definecolor{violetred}     {cmyk}{0   , 0.81, 0   , 0   }
\definecolor{rhodamine}     {cmyk}{0   , 0.82, 0   , 0   }
\definecolor{mulberry}      {cmyk}{0.34, 0.90, 0   , 0.02}
\definecolor{redviolet}     {cmyk}{0.07, 0.90, 0   , 0.34}
\definecolor{fuchsia}       {cmyk}{0.47, 0.91, 0   , 0.08}
\definecolor{lavender}      {cmyk}{0   , 0.48, 0   , 0   }
\definecolor{thistle}       {cmyk}{0.12, 0.59, 0   , 0   }
\definecolor{orchid}        {cmyk}{0.32, 0.64, 0   , 0   }
\definecolor{darkorchid}    {cmyk}{0.40, 0.80, 0.20, 0   }
\definecolor{purple}        {cmyk}{0.45, 0.86, 0   , 0   }
\definecolor{plum}          {cmyk}{0.50, 1   , 0   , 0   }
\definecolor{violet}        {cmyk}{0.79, 0.88, 0   , 0   }
\definecolor{royalpurple}   {cmyk}{0.75, 0.90, 0   , 0   }
\definecolor{blueviolet}    {cmyk}{0.86, 0.91, 0   , 0.04}
\definecolor{periwinkle}    {cmyk}{0.57, 0.55, 0   , 0   }
\definecolor{cadetblue}     {cmyk}{0.62, 0.57, 0.23, 0   }
\definecolor{cornflowerblue}{cmyk}{0.65, 0.13, 0   , 0   }
\definecolor{midnightblue}  {cmyk}{0.98, 0.13, 0   , 0.43}
\definecolor{navyblue}      {cmyk}{0.94, 0.54, 0   , 0   }
\definecolor{royalblue}     {cmyk}{1   , 0.50, 0   , 0   }
\definecolor{blue}          {cmyk}{1   , 1   , 0   , 0   }
\definecolor{cerulean}      {cmyk}{0.94, 0.11, 0   , 0   }
\definecolor{cyan}          {cmyk}{1   , 0   , 0   , 0   }
\definecolor{processblue}   {cmyk}{0.96, 0   , 0   , 0   }
\definecolor{skyblue}       {cmyk}{0.62, 0   , 0.12, 0   }
\definecolor{turquoise}     {cmyk}{0.85, 0   , 0.20, 0   }
\definecolor{tealblue}      {cmyk}{0.86, 0   , 0.34, 0.02}
\definecolor{aquamarine}    {cmyk}{0.82, 0   , 0.30, 0   }
\definecolor{bluegreen}     {cmyk}{0.85, 0   , 0.33, 0   }
\definecolor{emerald}       {cmyk}{1   , 0   , 0.50, 0   }
\definecolor{junglegreen}   {cmyk}{0.99, 0   , 0.52, 0   }
\definecolor{seagreen}      {cmyk}{0.69, 0   , 0.50, 0   }
\definecolor{green}         {cmyk}{1   , 0   , 1   , 0   }
\definecolor{forestgreen}   {cmyk}{0.91, 0   , 0.88, 0.12}
\definecolor{pinegreen}     {cmyk}{0.92, 0   , 0.59, 0.25}
\definecolor{limegreen}     {cmyk}{0.50, 0   , 1   , 0   }
\definecolor{yellowgreen}   {cmyk}{0.44, 0   , 0.74, 0   }
\definecolor{springgreen}   {cmyk}{0.26, 0   , 0.76, 0   }
\definecolor{olivegreen}    {cmyk}{0.64, 0   , 0.95, 0.40}
\definecolor{rawsienna}     {cmyk}{0   , 0.72, 1   , 0.45}
\definecolor{sepia}         {cmyk}{0   , 0.83, 1   , 0.70}
\definecolor{brown}         {cmyk}{0   , 0.81, 1   , 0.60}
\definecolor{tan}           {cmyk}{0.14, 0.42, 0.56, 0   }
\definecolor{gray}          {cmyk}{0   , 0   , 0   , 0.50}
\definecolor{black}         {cmyk}{0   , 0   , 0   , 1   }
\definecolor{white}         {cmyk}{0   , 0   , 0   , 0   } 

 \usepackage[]{hyperref} 
        
\usepackage{pgfplots}
\pgfplotsset{compat=newest}       
\usetikzlibrary{external}
 \usepackage{relinput}
 \usepackage{enumerate,graphicx}

\newtheorem{theorem}{Theorem}[section]

\newcounter{tikzsubfigcounter}[figure]
\renewcommand{\thetikzsubfigcounter}{\the\numexpr\value{figure}+1\relax\alph{tikzsubfigcounter}}

\newcounter{tikzsubfigcounterinvisible}[figure]
\renewcommand{\thetikzsubfigcounterinvisible}{\the\numexpr\value{figure}+1\relax\alph{tikzsubfigcounterinvisible}}

\newtheorem{thm}{\bf Theorem}[section]

\newtheorem{rem}[thm]{\bf  Remark}

\numberwithin{equation}{section}
\title{A numerical comparison of the method of moments for the population balance equation}
\author[lm]{Laura M\"uller}
\address[lm]{Fachbereich Mathematik, TU Kaiserslautern, Erwin-Schr\"odinger-Str., 67663 Kaiserslautern, Germany, {\tt lmueller@mathematik.uni-kl.de}}

\author[ak]{Axel Klar}
\address[ak]{Fachbereich Mathematik, TU Kaiserslautern, Erwin-Schr\"odinger-Str., 67663 Kaiserslautern, Germany, {\tt klar@mathematik.uni-kl.de}}

\author[fs]{Florian Schneider}
\address[fs]{Fachbereich Mathematik, TU Kaiserslautern, Erwin-Schr\"odinger-Str., 67663 Kaiserslautern, Germany, {\tt schneider@mathematik.uni-kl.de}}

\date{}
\input{colormap}
\input{usermacros}

\begin{document}

\begin{abstract}
We investigate the application of the method of moments approach for the one-dimensional population balance equation. 
We consider different types of moment closures, namely polynomial ($P_N$) closures,  maximum entropy ($M_N$) closures and the quadrature method of moments $QMOM_N$. Realizability issues and implementation details are discussed. 
The numerical examples range from spatially homogeneous cases to a population balance equation coupled with fluid dynamic equations for a lid-driven cavity test case. A detailed numerical discussion of accuracy, order of the moment method and computational time is given.
\end{abstract}
\maketitle

\noindent

% {\bf Key words.}

\input{Sections/introduction}
\input{Sections/PBE}
\input{Sections/momentmodels}

\input{Sections/modelling}

\input{Sections/results}

\input{Sections/outlook}

% Bibliography
%%%%%%%%%%%%%%
\bibliographystyle{siam}
\bibliography{bibliography}

\end{document}

%% file: colormap.tex
\pgfplotsset{colormap={mycolormap}{[1pt]
rgb(0pt)=(0.2081,0.1663,0.5292);          
rgb(1pt)=(0.208355,0.16778,0.532238);     
rgb(2pt)=(0.208611,0.169261,0.535275);    
rgb(3pt)=(0.208866,0.170741,0.538313);    
rgb(4pt)=(0.209121,0.172222,0.54135);     
rgb(5pt)=(0.209376,0.173702,0.544388);    
rgb(6pt)=(0.209632,0.175183,0.547425);    
rgb(7pt)=(0.209887,0.176663,0.550463);    
rgb(8pt)=(0.210134,0.178144,0.553505);    
rgb(9pt)=(0.210338,0.179624,0.556568);    
rgb(10pt)=(0.210542,0.181105,0.559631);   
rgb(11pt)=(0.210746,0.182585,0.562694);   
rgb(12pt)=(0.210944,0.184066,0.565763);   
rgb(13pt)=(0.211123,0.185546,0.568852);   
rgb(14pt)=(0.211302,0.187027,0.57194);    
rgb(15pt)=(0.21148,0.188507,0.575029);    
rgb(16pt)=(0.211642,0.189996,0.578117);   
rgb(17pt)=(0.21177,0.191502,0.581206);    
rgb(18pt)=(0.211897,0.193008,0.584295);   
rgb(19pt)=(0.212025,0.194514,0.587383);   
rgb(20pt)=(0.212132,0.19602,0.590472);    
rgb(21pt)=(0.212208,0.197526,0.59356);    
rgb(22pt)=(0.212285,0.199032,0.596649);   
rgb(23pt)=(0.212361,0.200538,0.599738);   
rgb(24pt)=(0.212413,0.202044,0.602839);   
rgb(25pt)=(0.212438,0.20355,0.605953);    
rgb(26pt)=(0.212464,0.205056,0.609067);   
rgb(27pt)=(0.212489,0.206562,0.612181);   
rgb(28pt)=(0.212471,0.208083,0.61531);    
rgb(29pt)=(0.21242,0.209614,0.61845);     
rgb(30pt)=(0.212368,0.211146,0.621589);   
rgb(31pt)=(0.212317,0.212677,0.624729);   
rgb(32pt)=(0.212216,0.214209,0.627868);   
rgb(33pt)=(0.212088,0.215741,0.631008);   
rgb(34pt)=(0.211961,0.217272,0.634148);   
rgb(35pt)=(0.211833,0.218804,0.637287);   
rgb(36pt)=(0.211668,0.220354,0.640446);   
rgb(37pt)=(0.211489,0.221911,0.643611);   
rgb(38pt)=(0.21131,0.223468,0.646776);    
rgb(39pt)=(0.211132,0.225025,0.649941);   
rgb(40pt)=(0.210848,0.226603,0.653107);   
rgb(41pt)=(0.210541,0.228186,0.656272);   
rgb(42pt)=(0.210235,0.229768,0.659437);   
rgb(43pt)=(0.209929,0.231351,0.662602);   
rgb(44pt)=(0.209553,0.232934,0.665767);   
rgb(45pt)=(0.20917,0.234516,0.668932);    
rgb(46pt)=(0.208787,0.236099,0.672098);   
rgb(47pt)=(0.208405,0.237681,0.675263);   
rgb(48pt)=(0.20787,0.239289,0.678453);    
rgb(49pt)=(0.207334,0.240897,0.681644);   
rgb(50pt)=(0.206798,0.242505,0.684835);   
rgb(51pt)=(0.206255,0.244114,0.688025);   
rgb(52pt)=(0.205617,0.245722,0.691216);   
rgb(53pt)=(0.204979,0.24733,0.694407);    
rgb(54pt)=(0.204341,0.248938,0.697597);   
rgb(55pt)=(0.203675,0.250554,0.700792);   
rgb(56pt)=(0.202858,0.252213,0.704008);   
rgb(57pt)=(0.202041,0.253872,0.707224);   
rgb(58pt)=(0.201225,0.255531,0.710441);   
rgb(59pt)=(0.200372,0.257184,0.713657);   
rgb(60pt)=(0.199402,0.258818,0.716873);   
rgb(61pt)=(0.198432,0.260452,0.720089);   
rgb(62pt)=(0.197462,0.262085,0.723305);   
rgb(63pt)=(0.196419,0.263735,0.726522);   
rgb(64pt)=(0.195219,0.26542,0.729738);    
rgb(65pt)=(0.19402,0.267105,0.732954);    
rgb(66pt)=(0.19282,0.268789,0.73617);     
rgb(67pt)=(0.191549,0.270474,0.739386);   
rgb(68pt)=(0.19017,0.272159,0.742603);    
rgb(69pt)=(0.188792,0.273843,0.745819);   
rgb(70pt)=(0.187414,0.275528,0.749035);   
rgb(71pt)=(0.1859,0.277237,0.752264);     
rgb(72pt)=(0.184241,0.278973,0.755505);   
rgb(73pt)=(0.182581,0.280709,0.758747);   
rgb(74pt)=(0.180922,0.282444,0.761989);   
rgb(75pt)=(0.179133,0.284209,0.765245);   
rgb(76pt)=(0.177244,0.285996,0.768512);   
rgb(77pt)=(0.175356,0.287783,0.77178);    
rgb(78pt)=(0.173467,0.289569,0.775047);   
rgb(79pt)=(0.171363,0.291406,0.778314);   
rgb(80pt)=(0.169142,0.293269,0.781581);   
rgb(81pt)=(0.166922,0.295132,0.784849);   
rgb(82pt)=(0.164701,0.296996,0.788116);   
rgb(83pt)=(0.162238,0.298934,0.791365);   
rgb(84pt)=(0.159686,0.300899,0.794606);   
rgb(85pt)=(0.157133,0.302865,0.797848);   
rgb(86pt)=(0.15458,0.30483,0.80109);      
rgb(87pt)=(0.151738,0.306858,0.804352);   
rgb(88pt)=(0.148828,0.3089,0.80762);      
rgb(89pt)=(0.145918,0.310942,0.810887);   
rgb(90pt)=(0.143008,0.312984,0.814154);   
rgb(91pt)=(0.139687,0.31514,0.81733);     
rgb(92pt)=(0.136318,0.31731,0.820495);    
rgb(93pt)=(0.132949,0.319479,0.82366);    
rgb(94pt)=(0.129579,0.321649,0.826826);   
rgb(95pt)=(0.125811,0.323918,0.829841);   
rgb(96pt)=(0.122033,0.32619,0.832853);    
rgb(97pt)=(0.118256,0.328462,0.835865);   
rgb(98pt)=(0.114458,0.330737,0.838862);   
rgb(99pt)=(0.110349,0.333059,0.841619);   
rgb(100pt)=(0.106239,0.335382,0.844376);  
rgb(101pt)=(0.102129,0.337705,0.847132);  
rgb(102pt)=(0.0979874,0.340021,0.849835); 
rgb(103pt)=(0.093648,0.342292,0.852209);  
rgb(104pt)=(0.0893087,0.344564,0.854583); 
rgb(105pt)=(0.0849694,0.346836,0.856957); 
rgb(106pt)=(0.08063,0.349091,0.859234);   
rgb(107pt)=(0.0762907,0.351286,0.861174); 
rgb(108pt)=(0.0719514,0.353481,0.863114); 
rgb(109pt)=(0.067612,0.355676,0.865053);  
rgb(110pt)=(0.0633195,0.357817,0.866853); 
rgb(111pt)=(0.0591333,0.359833,0.868333); 
rgb(112pt)=(0.0549471,0.36185,0.869814);  
rgb(113pt)=(0.050761,0.363866,0.871294);  
rgb(114pt)=(0.0466838,0.365823,0.872626); 
rgb(115pt)=(0.0427784,0.367687,0.873724); 
rgb(116pt)=(0.038873,0.36955,0.874821);   
rgb(117pt)=(0.0349676,0.371414,0.875919); 
rgb(118pt)=(0.0315066,0.373217,0.876872); 
rgb(119pt)=(0.0285456,0.374953,0.877664); 
rgb(120pt)=(0.0255847,0.376688,0.878455); 
rgb(121pt)=(0.0226237,0.378424,0.879246); 
rgb(122pt)=(0.0202132,0.380061,0.879868); 
rgb(123pt)=(0.0182477,0.381618,0.880353); 
rgb(124pt)=(0.0162823,0.383175,0.880838); 
rgb(125pt)=(0.0143168,0.384732,0.881323); 
rgb(126pt)=(0.0127892,0.386241,0.881695); 
rgb(127pt)=(0.0115129,0.387721,0.882001); 
rgb(128pt)=(0.0102366,0.389202,0.882307); 
rgb(129pt)=(0.00896036,0.390682,0.882614);
rgb(130pt)=(0.00812372,0.392089,0.88281); 
rgb(131pt)=(0.00746006,0.393468,0.882963);
rgb(132pt)=(0.0067964,0.394846,0.883116); 
rgb(133pt)=(0.00613273,0.396224,0.883269);
rgb(134pt)=(0.00581622,0.397562,0.88332); 
rgb(135pt)=(0.00558649,0.398889,0.883346);
rgb(136pt)=(0.00535676,0.400217,0.883371);
rgb(137pt)=(0.00512703,0.401544,0.883397);
rgb(138pt)=(0.00516757,0.402804,0.883332);
rgb(139pt)=(0.00524414,0.404054,0.883256);
rgb(140pt)=(0.00532072,0.405305,0.883179);
rgb(141pt)=(0.0053973,0.406556,0.883103); 
rgb(142pt)=(0.00572012,0.407757,0.882952);
rgb(143pt)=(0.00605195,0.408957,0.882799);
rgb(144pt)=(0.00638378,0.410157,0.882646);
rgb(145pt)=(0.00672643,0.411355,0.882489);
rgb(146pt)=(0.00728799,0.412529,0.882259);
rgb(147pt)=(0.00784955,0.413704,0.88203); 
rgb(148pt)=(0.00841111,0.414878,0.8818);  
rgb(149pt)=(0.00898919,0.416045,0.881564);
rgb(150pt)=(0.00967838,0.417168,0.881283);
rgb(151pt)=(0.0103676,0.418292,0.881002); 
rgb(152pt)=(0.0110568,0.419415,0.880721); 
rgb(153pt)=(0.011773,0.420532,0.880435);  
rgb(154pt)=(0.0125898,0.42163,0.880129);  
rgb(155pt)=(0.0134066,0.422728,0.879823); 
rgb(156pt)=(0.0142234,0.423825,0.879516); 
rgb(157pt)=(0.0150703,0.424915,0.879195); 
rgb(158pt)=(0.0159892,0.425987,0.878838); 
rgb(159pt)=(0.0169081,0.427059,0.87848);  
rgb(160pt)=(0.017827,0.428132,0.878123);  
rgb(161pt)=(0.0187748,0.429194,0.877746); 
rgb(162pt)=(0.0197703,0.430241,0.877338); 
rgb(163pt)=(0.0207658,0.431287,0.876929); 
rgb(164pt)=(0.0217613,0.432334,0.876521); 
rgb(165pt)=(0.0227802,0.43338,0.876113);  
rgb(166pt)=(0.0238267,0.434427,0.875704); 
rgb(167pt)=(0.0248733,0.435473,0.875296); 
rgb(168pt)=(0.0259198,0.43652,0.874887);  
rgb(169pt)=(0.0269802,0.437553,0.874451); 
rgb(170pt)=(0.0280523,0.438574,0.873992); 
rgb(171pt)=(0.0291243,0.439595,0.873532); 
rgb(172pt)=(0.0301964,0.440616,0.873073); 
rgb(173pt)=(0.0312844,0.441621,0.872614); 
rgb(174pt)=(0.032382,0.442616,0.872154);  
rgb(175pt)=(0.0334796,0.443612,0.871695); 
rgb(176pt)=(0.0345772,0.444607,0.871235); 
rgb(177pt)=(0.0357108,0.445603,0.870758); 
rgb(178pt)=(0.0368595,0.446598,0.870273); 
rgb(179pt)=(0.0380081,0.447594,0.869788); 
rgb(180pt)=(0.0391568,0.448589,0.869303); 
rgb(181pt)=(0.0402652,0.449565,0.868798); 
rgb(182pt)=(0.0413628,0.450535,0.868287); 
rgb(183pt)=(0.0424604,0.451505,0.867777); 
rgb(184pt)=(0.043558,0.452474,0.867266);  
rgb(185pt)=(0.0445889,0.453444,0.866756); 
rgb(186pt)=(0.0456099,0.454414,0.866245); 
rgb(187pt)=(0.0466309,0.455384,0.865735); 
rgb(188pt)=(0.047652,0.456354,0.865224);  
rgb(189pt)=(0.0486,0.457324,0.864714);    
rgb(190pt)=(0.0495444,0.458294,0.864203); 
rgb(191pt)=(0.0504889,0.459264,0.863692); 
rgb(192pt)=(0.0514315,0.460234,0.863181); 
rgb(193pt)=(0.0523249,0.461204,0.862645); 
rgb(194pt)=(0.0532183,0.462174,0.862109); 
rgb(195pt)=(0.0541117,0.463144,0.861573); 
rgb(196pt)=(0.0549991,0.464111,0.861034); 
rgb(197pt)=(0.0558414,0.465056,0.860472); 
rgb(198pt)=(0.0566838,0.466,0.859911);    
rgb(199pt)=(0.0575261,0.466944,0.859349); 
rgb(200pt)=(0.0583532,0.467889,0.858793); 
rgb(201pt)=(0.0591189,0.468833,0.858257); 
rgb(202pt)=(0.0598847,0.469778,0.857721); 
rgb(203pt)=(0.0606505,0.470722,0.857185); 
rgb(204pt)=(0.0614018,0.471667,0.856641); 
rgb(205pt)=(0.0621165,0.472611,0.85608);  
rgb(206pt)=(0.0628312,0.473556,0.855518); 
rgb(207pt)=(0.0635459,0.4745,0.854957);   
rgb(208pt)=(0.064242,0.475444,0.854405);  
rgb(209pt)=(0.0649057,0.476389,0.853868); 
rgb(210pt)=(0.0655694,0.477333,0.853332); 
rgb(211pt)=(0.066233,0.478278,0.852796);  
rgb(212pt)=(0.0668625,0.479222,0.852249); 
rgb(213pt)=(0.0674495,0.480167,0.851687); 
rgb(214pt)=(0.0680366,0.481111,0.851126); 
rgb(215pt)=(0.0686237,0.482056,0.850564); 
rgb(216pt)=(0.0691838,0.483,0.850003);    
rgb(217pt)=(0.0697198,0.483944,0.849441); 
rgb(218pt)=(0.0702559,0.484889,0.84888);  
rgb(219pt)=(0.0707919,0.485833,0.848318); 
rgb(220pt)=(0.0712967,0.486778,0.847772); 
rgb(221pt)=(0.0717817,0.487722,0.847236); 
rgb(222pt)=(0.0722667,0.488667,0.8467);   
rgb(223pt)=(0.0727517,0.489611,0.846164); 
rgb(224pt)=(0.0732012,0.490573,0.845628); 
rgb(225pt)=(0.0736351,0.491543,0.845092); 
rgb(226pt)=(0.0740691,0.492513,0.844556); 
rgb(227pt)=(0.074503,0.493483,0.84402);   
rgb(228pt)=(0.0748973,0.494433,0.843484); 
rgb(229pt)=(0.0752802,0.495378,0.842948); 
rgb(230pt)=(0.0756631,0.496322,0.842412); 
rgb(231pt)=(0.0760459,0.497267,0.841876); 
rgb(232pt)=(0.0763631,0.498233,0.841362); 
rgb(233pt)=(0.0766694,0.499203,0.840851); 
rgb(234pt)=(0.0769757,0.500173,0.840341); 
rgb(235pt)=(0.077282,0.501143,0.83983);   
rgb(236pt)=(0.0775162,0.502137,0.83932);  
rgb(237pt)=(0.0777459,0.503132,0.838809); 
rgb(238pt)=(0.0779757,0.504128,0.838298); 
rgb(239pt)=(0.0782042,0.505123,0.837789); 
rgb(240pt)=(0.0783829,0.506093,0.837304); 
rgb(241pt)=(0.0785616,0.507063,0.836819); 
rgb(242pt)=(0.0787402,0.508033,0.836334); 
rgb(243pt)=(0.0789135,0.509008,0.835851); 
rgb(244pt)=(0.0790411,0.510029,0.835392); 
rgb(245pt)=(0.0791688,0.51105,0.834932);  
rgb(246pt)=(0.0792964,0.512071,0.834473); 
rgb(247pt)=(0.0794048,0.513092,0.834018); 
rgb(248pt)=(0.0794303,0.514113,0.833584); 
rgb(249pt)=(0.0794559,0.515134,0.83315);  
rgb(250pt)=(0.0794814,0.516155,0.832717); 
rgb(251pt)=(0.0794862,0.517183,0.832289); 
rgb(252pt)=(0.0794351,0.51823,0.831881);  
rgb(253pt)=(0.0793841,0.519276,0.831473); 
rgb(254pt)=(0.079333,0.520323,0.831064);  
rgb(255pt)=(0.079255,0.521369,0.830665);  
rgb(256pt)=(0.0791273,0.522416,0.830282); 
rgb(257pt)=(0.0789997,0.523462,0.829899); 
rgb(258pt)=(0.0788721,0.524509,0.829516); 
rgb(259pt)=(0.0786889,0.525589,0.829156); 
rgb(260pt)=(0.0784336,0.526712,0.828824); 
rgb(261pt)=(0.0781784,0.527835,0.828492); 
rgb(262pt)=(0.0779231,0.528958,0.82816);  
rgb(263pt)=(0.077615,0.530081,0.827868);  
rgb(264pt)=(0.0772577,0.531205,0.827613); 
rgb(265pt)=(0.0769003,0.532328,0.827357); 
rgb(266pt)=(0.0765429,0.533451,0.827102); 
rgb(267pt)=(0.0761243,0.534589,0.826862); 
rgb(268pt)=(0.0756649,0.535738,0.826632); 
rgb(269pt)=(0.0752054,0.536886,0.826403); 
rgb(270pt)=(0.0747459,0.538035,0.826173); 
rgb(271pt)=(0.0742168,0.539219,0.825961); 
rgb(272pt)=(0.0736553,0.540418,0.825756); 
rgb(273pt)=(0.0730937,0.541618,0.825552); 
rgb(274pt)=(0.0725321,0.542818,0.825348); 
rgb(275pt)=(0.0718925,0.544037,0.825183); 
rgb(276pt)=(0.0712288,0.545262,0.82503);  
rgb(277pt)=(0.0705652,0.546487,0.824877); 
rgb(278pt)=(0.0699015,0.547713,0.824723); 
rgb(279pt)=(0.0691514,0.548938,0.824614); 
rgb(280pt)=(0.0683856,0.550163,0.824511); 
rgb(281pt)=(0.0676198,0.551388,0.824409); 
rgb(282pt)=(0.0668541,0.552614,0.824307); 
rgb(283pt)=(0.0660408,0.553886,0.824205); 
rgb(284pt)=(0.065224,0.555162,0.824103);  
rgb(285pt)=(0.0644072,0.556439,0.824001); 
rgb(286pt)=(0.0635892,0.557715,0.823899); 
rgb(287pt)=(0.0626703,0.558991,0.823848); 
rgb(288pt)=(0.0617514,0.560268,0.823797); 
rgb(289pt)=(0.0608324,0.561544,0.823746); 
rgb(290pt)=(0.0599087,0.56282,0.823693);  
rgb(291pt)=(0.0589387,0.564096,0.823616); 
rgb(292pt)=(0.0579688,0.565373,0.82354);  
rgb(293pt)=(0.0569988,0.566649,0.823463); 
rgb(294pt)=(0.0560243,0.567925,0.823386); 
rgb(295pt)=(0.0550288,0.569202,0.82331);  
rgb(296pt)=(0.0540333,0.570478,0.823233); 
rgb(297pt)=(0.0530378,0.571754,0.823157); 
rgb(298pt)=(0.0520423,0.57303,0.82308);   
rgb(299pt)=(0.0510468,0.574307,0.823004); 
rgb(300pt)=(0.0500514,0.575583,0.822927); 
rgb(301pt)=(0.0490559,0.576859,0.82285);  
rgb(302pt)=(0.0480604,0.578127,0.822756); 
rgb(303pt)=(0.0470649,0.579377,0.822629); 
rgb(304pt)=(0.0460694,0.580628,0.822501); 
rgb(305pt)=(0.0450739,0.581879,0.822374); 
rgb(306pt)=(0.0441,0.583119,0.822235);    
rgb(307pt)=(0.0431556,0.584344,0.822082); 
rgb(308pt)=(0.0422111,0.585569,0.821929); 
rgb(309pt)=(0.0412667,0.586795,0.821776); 
rgb(310pt)=(0.0403351,0.58802,0.821597);  
rgb(311pt)=(0.0394162,0.589245,0.821392); 
rgb(312pt)=(0.0384973,0.59047,0.821188);  
rgb(313pt)=(0.0375784,0.591695,0.820984); 
rgb(314pt)=(0.0367495,0.592891,0.820735); 
rgb(315pt)=(0.0359838,0.594065,0.820454); 
rgb(316pt)=(0.035218,0.595239,0.820173);  
rgb(317pt)=(0.0344523,0.596413,0.819892); 
rgb(318pt)=(0.0337721,0.597553,0.819595); 
rgb(319pt)=(0.0331339,0.598676,0.819288); 
rgb(320pt)=(0.0324958,0.599799,0.818982); 
rgb(321pt)=(0.0318577,0.600923,0.818676); 
rgb(322pt)=(0.0312964,0.602026,0.818312); 
rgb(323pt)=(0.0307604,0.603124,0.817929); 
rgb(324pt)=(0.0302243,0.604222,0.817546); 
rgb(325pt)=(0.0296883,0.605319,0.817163); 
rgb(326pt)=(0.0292375,0.606395,0.816738); 
rgb(327pt)=(0.0288036,0.607468,0.816304); 
rgb(328pt)=(0.0283697,0.60854,0.81587);   
rgb(329pt)=(0.0279357,0.609612,0.815436); 
rgb(330pt)=(0.0275721,0.610637,0.814955); 
rgb(331pt)=(0.0272147,0.611658,0.81447);  
rgb(332pt)=(0.0268574,0.612679,0.813985); 
rgb(333pt)=(0.0265,0.6137,0.8135);        
rgb(334pt)=(0.0262447,0.614695,0.812964); 
rgb(335pt)=(0.0259895,0.615691,0.812428); 
rgb(336pt)=(0.0257342,0.616686,0.811892); 
rgb(337pt)=(0.0254853,0.61768,0.811352);  
rgb(338pt)=(0.0253066,0.61865,0.810765);  
rgb(339pt)=(0.0251279,0.61962,0.810177);  
rgb(340pt)=(0.0249492,0.62059,0.80959);   
rgb(341pt)=(0.024779,0.621551,0.808995);  
rgb(342pt)=(0.0246514,0.62247,0.808357);  
rgb(343pt)=(0.0245237,0.623389,0.807719); 
rgb(344pt)=(0.0243961,0.624308,0.80708);  
rgb(345pt)=(0.0242748,0.625221,0.80643);  
rgb(346pt)=(0.0241727,0.626114,0.805741); 
rgb(347pt)=(0.0240706,0.627008,0.805051); 
rgb(348pt)=(0.0239685,0.627901,0.804362); 
rgb(349pt)=(0.0238832,0.628786,0.803656); 
rgb(350pt)=(0.0238321,0.629654,0.802916); 
rgb(351pt)=(0.0237811,0.630522,0.802176); 
rgb(352pt)=(0.02373,0.631389,0.801435);   
rgb(353pt)=(0.023679,0.632247,0.800685);  
rgb(354pt)=(0.0236279,0.633089,0.799919); 
rgb(355pt)=(0.0235769,0.633932,0.799153); 
rgb(356pt)=(0.0235258,0.634774,0.798387); 
rgb(357pt)=(0.0234748,0.635604,0.797596); 
rgb(358pt)=(0.0234237,0.63642,0.79678);   
rgb(359pt)=(0.0233727,0.637237,0.795963); 
rgb(360pt)=(0.0233216,0.638054,0.795146); 
rgb(361pt)=(0.0232706,0.638856,0.794329); 
rgb(362pt)=(0.0232195,0.639647,0.793512); 
rgb(363pt)=(0.0231685,0.640439,0.792695); 
rgb(364pt)=(0.0231174,0.64123,0.791879);  
rgb(365pt)=(0.0230832,0.642005,0.791011); 
rgb(366pt)=(0.0230577,0.64277,0.790118);  
rgb(367pt)=(0.0230321,0.643536,0.789225); 
rgb(368pt)=(0.0230066,0.644302,0.788331); 
rgb(369pt)=(0.0229811,0.645049,0.787438); 
rgb(370pt)=(0.0229556,0.645789,0.786544); 
rgb(371pt)=(0.02293,0.646529,0.785651);   
rgb(372pt)=(0.0229045,0.647269,0.784758); 
rgb(373pt)=(0.022858,0.64801,0.783843);   
rgb(374pt)=(0.0228069,0.64875,0.782924);  
rgb(375pt)=(0.0227559,0.64949,0.782005);  
rgb(376pt)=(0.0227048,0.65023,0.781086);  
rgb(377pt)=(0.0227,0.650947,0.780144);    
rgb(378pt)=(0.0227,0.651662,0.7792);      
rgb(379pt)=(0.0227,0.652377,0.778256);    
rgb(380pt)=(0.0227,0.653092,0.777311);    
rgb(381pt)=(0.0228261,0.653781,0.776341); 
rgb(382pt)=(0.0229538,0.65447,0.775371);  
rgb(383pt)=(0.0230814,0.655159,0.774402); 
rgb(384pt)=(0.0232108,0.655849,0.77343);  
rgb(385pt)=(0.023364,0.656538,0.772434);  
rgb(386pt)=(0.0235171,0.657227,0.771439); 
rgb(387pt)=(0.0236703,0.657916,0.770443); 
rgb(388pt)=(0.0238312,0.658602,0.769444); 
rgb(389pt)=(0.0240354,0.659265,0.768423); 
rgb(390pt)=(0.0242396,0.659929,0.767402); 
rgb(391pt)=(0.0244438,0.660592,0.766381); 
rgb(392pt)=(0.0247021,0.661256,0.765354); 
rgb(393pt)=(0.025136,0.66192,0.764307);   
rgb(394pt)=(0.02557,0.662583,0.763261);   
rgb(395pt)=(0.0260039,0.663247,0.762214); 
rgb(396pt)=(0.0264541,0.663911,0.761168); 
rgb(397pt)=(0.026939,0.664574,0.760121);  
rgb(398pt)=(0.027424,0.665238,0.759074);  
rgb(399pt)=(0.027909,0.665902,0.758028);  
rgb(400pt)=(0.028445,0.666555,0.756971);  
rgb(401pt)=(0.0290577,0.667193,0.755899); 
rgb(402pt)=(0.0296703,0.667832,0.754827); 
rgb(403pt)=(0.0302829,0.66847,0.753755);  
rgb(404pt)=(0.030994,0.669095,0.752683);  
rgb(405pt)=(0.0318108,0.669708,0.751611); 
rgb(406pt)=(0.0326276,0.670321,0.750539); 
rgb(407pt)=(0.0334444,0.670933,0.749467); 
rgb(408pt)=(0.0343045,0.67156,0.748366);  
rgb(409pt)=(0.0351979,0.672198,0.747243); 
rgb(410pt)=(0.0360913,0.672837,0.74612);  
rgb(411pt)=(0.0369847,0.673475,0.744996); 
rgb(412pt)=(0.0380432,0.674096,0.743873); 
rgb(413pt)=(0.0391919,0.674709,0.74275);  
rgb(414pt)=(0.0403405,0.675322,0.741627); 
rgb(415pt)=(0.0414892,0.675934,0.740504); 
rgb(416pt)=(0.0427123,0.676528,0.739381); 
rgb(417pt)=(0.0439631,0.677115,0.738258); 
rgb(418pt)=(0.0452138,0.677702,0.737135); 
rgb(419pt)=(0.0464646,0.678289,0.736011); 
rgb(420pt)=(0.0477153,0.678897,0.734868); 
rgb(421pt)=(0.0489661,0.67951,0.733719);  
rgb(422pt)=(0.0502168,0.680123,0.73257);  
rgb(423pt)=(0.0514676,0.680735,0.731422); 
rgb(424pt)=(0.0529237,0.681325,0.73025);  
rgb(425pt)=(0.0544042,0.681912,0.729076); 
rgb(426pt)=(0.0558847,0.682499,0.727902); 
rgb(427pt)=(0.0573652,0.683086,0.726728); 
rgb(428pt)=(0.0587709,0.683673,0.725553); 
rgb(429pt)=(0.0601748,0.68426,0.724379);  
rgb(430pt)=(0.0615787,0.684847,0.723205); 
rgb(431pt)=(0.0629946,0.685435,0.722028); 
rgb(432pt)=(0.0646027,0.686022,0.720803); 
rgb(433pt)=(0.0662108,0.686609,0.719577); 
rgb(434pt)=(0.0678189,0.687196,0.718352); 
rgb(435pt)=(0.069427,0.687779,0.717131);  
rgb(436pt)=(0.0710351,0.688341,0.715931); 
rgb(437pt)=(0.0726432,0.688902,0.714731); 
rgb(438pt)=(0.0742514,0.689464,0.713532); 
rgb(439pt)=(0.0758709,0.690026,0.712326); 
rgb(440pt)=(0.07753,0.690587,0.711101);   
rgb(441pt)=(0.0791892,0.691149,0.709876); 
rgb(442pt)=(0.0808483,0.69171,0.70865);   
rgb(443pt)=(0.0825387,0.692272,0.707417); 
rgb(444pt)=(0.0843,0.692833,0.706167);    
rgb(445pt)=(0.0860613,0.693395,0.704916); 
rgb(446pt)=(0.0878225,0.693956,0.703665); 
rgb(447pt)=(0.089564,0.694518,0.702405);  
rgb(448pt)=(0.0912742,0.69508,0.701128);  
rgb(449pt)=(0.0929844,0.695641,0.699852); 
rgb(450pt)=(0.0946946,0.696203,0.698576); 
rgb(451pt)=(0.0965009,0.696752,0.697299); 
rgb(452pt)=(0.0984153,0.697288,0.696023); 
rgb(453pt)=(0.10033,0.697824,0.694747);   
rgb(454pt)=(0.102244,0.69836,0.693471);   
rgb(455pt)=(0.10413,0.698896,0.69218);    
rgb(456pt)=(0.105994,0.699432,0.690878);  
rgb(457pt)=(0.107857,0.699968,0.689577);  
rgb(458pt)=(0.10972,0.700505,0.688275);   
rgb(459pt)=(0.111632,0.701041,0.686973);  
rgb(460pt)=(0.113572,0.701577,0.685671);  
rgb(461pt)=(0.115512,0.702113,0.684369);  
rgb(462pt)=(0.117452,0.702649,0.683068);  
rgb(463pt)=(0.119429,0.703185,0.681747);  
rgb(464pt)=(0.12142,0.703721,0.68042);    
rgb(465pt)=(0.123411,0.704257,0.679093);  
rgb(466pt)=(0.125402,0.704793,0.677765);  
rgb(467pt)=(0.127372,0.705308,0.676438);  
rgb(468pt)=(0.129338,0.705819,0.675111);  
rgb(469pt)=(0.131303,0.706329,0.673783);  
rgb(470pt)=(0.133269,0.70684,0.672456);   
rgb(471pt)=(0.135369,0.70735,0.671084);   
rgb(472pt)=(0.137488,0.707861,0.669705);  
rgb(473pt)=(0.139607,0.708371,0.668327);  
rgb(474pt)=(0.141725,0.708882,0.666949);  
rgb(475pt)=(0.143795,0.709392,0.665595);  
rgb(476pt)=(0.145862,0.709903,0.664242);  
rgb(477pt)=(0.14793,0.710414,0.662889);   
rgb(478pt)=(0.150003,0.710924,0.661534);  
rgb(479pt)=(0.152198,0.711435,0.66013);   
rgb(480pt)=(0.154394,0.711945,0.658726);  
rgb(481pt)=(0.156589,0.712456,0.657322);  
rgb(482pt)=(0.158784,0.712963,0.655922);  
rgb(483pt)=(0.160979,0.713448,0.654543);  
rgb(484pt)=(0.163174,0.713933,0.653165);  
rgb(485pt)=(0.16537,0.714418,0.651786);   
rgb(486pt)=(0.16757,0.714908,0.650397);   
rgb(487pt)=(0.169791,0.715419,0.648968);  
rgb(488pt)=(0.172012,0.715929,0.647538);  
rgb(489pt)=(0.174232,0.71644,0.646109);   
rgb(490pt)=(0.176483,0.716935,0.64468);   
rgb(491pt)=(0.178806,0.717395,0.64325);   
rgb(492pt)=(0.181129,0.717854,0.641821);  
rgb(493pt)=(0.183452,0.718314,0.640391);  
rgb(494pt)=(0.185755,0.718783,0.638952);  
rgb(495pt)=(0.188027,0.719268,0.637497);  
rgb(496pt)=(0.190299,0.719753,0.636042);  
rgb(497pt)=(0.192571,0.720238,0.634587);  
rgb(498pt)=(0.194913,0.720711,0.633132);  
rgb(499pt)=(0.197338,0.72117,0.631677);   
rgb(500pt)=(0.199762,0.72163,0.630223);   
rgb(501pt)=(0.202187,0.722089,0.628768);  
rgb(502pt)=(0.204612,0.722549,0.627299);  
rgb(503pt)=(0.207037,0.723008,0.625818);  
rgb(504pt)=(0.209462,0.723468,0.624338);  
rgb(505pt)=(0.211887,0.723927,0.622857);  
rgb(506pt)=(0.214328,0.724386,0.621377);  
rgb(507pt)=(0.216778,0.724846,0.619896);  
rgb(508pt)=(0.219229,0.725305,0.618416);  
rgb(509pt)=(0.221679,0.725765,0.616935);  
rgb(510pt)=(0.224202,0.726188,0.615455);  
rgb(511pt)=(0.226754,0.726597,0.613974);  
rgb(512pt)=(0.229307,0.727005,0.612494);  
rgb(513pt)=(0.231859,0.727414,0.611014);  
rgb(514pt)=(0.234392,0.727842,0.609513);  
rgb(515pt)=(0.236919,0.728276,0.608007);  
rgb(516pt)=(0.239446,0.72871,0.606501);   
rgb(517pt)=(0.241973,0.729144,0.604995);  
rgb(518pt)=(0.244611,0.729556,0.603467);  
rgb(519pt)=(0.247266,0.729964,0.601935);  
rgb(520pt)=(0.24992,0.730372,0.600404);   
rgb(521pt)=(0.252575,0.730781,0.598872);  
rgb(522pt)=(0.25523,0.731189,0.597365);   
rgb(523pt)=(0.257884,0.731598,0.595859);  
rgb(524pt)=(0.260539,0.732006,0.594353);  
rgb(525pt)=(0.263194,0.732414,0.592846);  
rgb(526pt)=(0.265848,0.732796,0.591314);  
rgb(527pt)=(0.268503,0.733179,0.589783);  
rgb(528pt)=(0.271158,0.733562,0.588251);  
rgb(529pt)=(0.27383,0.733945,0.58672);    
rgb(530pt)=(0.276638,0.734328,0.585188);  
rgb(531pt)=(0.279446,0.734711,0.583657);  
rgb(532pt)=(0.282254,0.735094,0.582125);  
rgb(533pt)=(0.285051,0.735471,0.580594);  
rgb(534pt)=(0.287808,0.735829,0.579062);  
rgb(535pt)=(0.290565,0.736186,0.577531);  
rgb(536pt)=(0.293322,0.736544,0.575999);  
rgb(537pt)=(0.2961,0.736894,0.574468);    
rgb(538pt)=(0.298933,0.737226,0.572936);  
rgb(539pt)=(0.301767,0.737557,0.571405);  
rgb(540pt)=(0.3046,0.737889,0.569873);    
rgb(541pt)=(0.307452,0.738221,0.568351);  
rgb(542pt)=(0.310336,0.738553,0.566845);  
rgb(543pt)=(0.313221,0.738885,0.565339);  
rgb(544pt)=(0.316105,0.739217,0.563833);  
rgb(545pt)=(0.318978,0.739537,0.562315);  
rgb(546pt)=(0.321837,0.739843,0.560784);  
rgb(547pt)=(0.324696,0.74015,0.559252);   
rgb(548pt)=(0.327555,0.740456,0.557721);  
rgb(549pt)=(0.330468,0.740749,0.556216);  
rgb(550pt)=(0.333429,0.741029,0.554736);  
rgb(551pt)=(0.336389,0.74131,0.553255);   
rgb(552pt)=(0.33935,0.741591,0.551775);   
rgb(553pt)=(0.342296,0.741872,0.550279);  
rgb(554pt)=(0.345231,0.742153,0.548773);  
rgb(555pt)=(0.348167,0.742433,0.547267);  
rgb(556pt)=(0.351102,0.742714,0.545761);  
rgb(557pt)=(0.354038,0.742977,0.54429);   
rgb(558pt)=(0.356973,0.743232,0.542835);  
rgb(559pt)=(0.359908,0.743488,0.54138);   
rgb(560pt)=(0.362844,0.743743,0.539925);  
rgb(561pt)=(0.365839,0.743959,0.53847);   
rgb(562pt)=(0.368851,0.744163,0.537015);  
rgb(563pt)=(0.371863,0.744367,0.53556);   
rgb(564pt)=(0.374875,0.744571,0.534105);  
rgb(565pt)=(0.377843,0.744775,0.532672);  
rgb(566pt)=(0.380804,0.74498,0.531243);   
rgb(567pt)=(0.383765,0.745184,0.529814);  
rgb(568pt)=(0.386726,0.745388,0.528384);  
rgb(569pt)=(0.389711,0.745568,0.527003);  
rgb(570pt)=(0.392697,0.745747,0.525624);  
rgb(571pt)=(0.395684,0.745926,0.524246);  
rgb(572pt)=(0.39867,0.746104,0.522868);   
rgb(573pt)=(0.401657,0.746257,0.521489);  
rgb(574pt)=(0.404643,0.74641,0.520111);   
rgb(575pt)=(0.40763,0.746563,0.518732);   
rgb(576pt)=(0.410611,0.746716,0.517359);  
rgb(577pt)=(0.413546,0.746869,0.516032);  
rgb(578pt)=(0.416482,0.747023,0.514705);  
rgb(579pt)=(0.419417,0.747176,0.513377);  
rgb(580pt)=(0.422357,0.747319,0.512055);  
rgb(581pt)=(0.425318,0.747421,0.510753);  
rgb(582pt)=(0.428279,0.747523,0.509451);  
rgb(583pt)=(0.43124,0.747626,0.50815);    
rgb(584pt)=(0.43418,0.747735,0.506848);   
rgb(585pt)=(0.437065,0.747862,0.505546);  
rgb(586pt)=(0.439949,0.74799,0.504244);   
rgb(587pt)=(0.442834,0.748117,0.502942);  
rgb(588pt)=(0.445727,0.748227,0.501659);  
rgb(589pt)=(0.448637,0.748304,0.500408);  
rgb(590pt)=(0.451547,0.74838,0.499157);   
rgb(591pt)=(0.454457,0.748457,0.497906);  
rgb(592pt)=(0.457333,0.748522,0.496667);  
rgb(593pt)=(0.460167,0.748573,0.495441);  
rgb(594pt)=(0.463,0.748624,0.494216);     
rgb(595pt)=(0.465833,0.748675,0.492991);  
rgb(596pt)=(0.468667,0.748726,0.491779);  
rgb(597pt)=(0.4715,0.748777,0.490579);    
rgb(598pt)=(0.474333,0.748829,0.48938);   
rgb(599pt)=(0.477167,0.74888,0.48818);    
rgb(600pt)=(0.479969,0.748931,0.486995);  
rgb(601pt)=(0.482752,0.748982,0.485821);  
rgb(602pt)=(0.485534,0.749033,0.484647);  
rgb(603pt)=(0.488316,0.749084,0.483473);  
rgb(604pt)=(0.491081,0.7491,0.482316);    
rgb(605pt)=(0.493838,0.7491,0.481168);    
rgb(606pt)=(0.496595,0.7491,0.480019);    
rgb(607pt)=(0.499351,0.7491,0.47887);     
rgb(608pt)=(0.502069,0.74912,0.477722);   
rgb(609pt)=(0.504775,0.749145,0.476573);  
rgb(610pt)=(0.50748,0.749171,0.475424);   
rgb(611pt)=(0.510186,0.749196,0.474276);  
rgb(612pt)=(0.512892,0.7492,0.47317);     
rgb(613pt)=(0.515598,0.7492,0.472073);    
rgb(614pt)=(0.518303,0.7492,0.470975);    
rgb(615pt)=(0.521009,0.7492,0.469877);    
rgb(616pt)=(0.523644,0.749176,0.46878);   
rgb(617pt)=(0.526273,0.749151,0.467682);  
rgb(618pt)=(0.528902,0.749125,0.466585);  
rgb(619pt)=(0.531531,0.7491,0.465487);    
rgb(620pt)=(0.53416,0.749074,0.464415);   
rgb(621pt)=(0.536789,0.749049,0.463343);  
rgb(622pt)=(0.539418,0.749023,0.462271);  
rgb(623pt)=(0.542043,0.748998,0.461199);  
rgb(624pt)=(0.544621,0.748972,0.460127);  
rgb(625pt)=(0.547199,0.748947,0.459055);  
rgb(626pt)=(0.549777,0.748921,0.457983);  
rgb(627pt)=(0.55235,0.748891,0.45692);    
rgb(628pt)=(0.554903,0.74884,0.455899);   
rgb(629pt)=(0.557456,0.748789,0.454878);  
rgb(630pt)=(0.560008,0.748738,0.453857);  
rgb(631pt)=(0.562554,0.748687,0.452829);  
rgb(632pt)=(0.565081,0.748636,0.451783);  
rgb(633pt)=(0.567608,0.748585,0.450736);  
rgb(634pt)=(0.570135,0.748534,0.449689);  
rgb(635pt)=(0.572653,0.748474,0.44866);   
rgb(636pt)=(0.575155,0.748397,0.447665);  
rgb(637pt)=(0.577656,0.748321,0.446669);  
rgb(638pt)=(0.580158,0.748244,0.445674);  
rgb(639pt)=(0.582649,0.748168,0.444678);  
rgb(640pt)=(0.585125,0.748091,0.443683);  
rgb(641pt)=(0.587601,0.748014,0.442687);  
rgb(642pt)=(0.590077,0.747938,0.441692);  
rgb(643pt)=(0.59254,0.747861,0.440709);   
rgb(644pt)=(0.59499,0.747785,0.439739);   
rgb(645pt)=(0.597441,0.747708,0.438769);  
rgb(646pt)=(0.599891,0.747632,0.437799);  
rgb(647pt)=(0.602311,0.747555,0.436814);  
rgb(648pt)=(0.604711,0.747478,0.435819);  
rgb(649pt)=(0.60711,0.747402,0.434823);   
rgb(650pt)=(0.60951,0.747325,0.433828);   
rgb(651pt)=(0.611909,0.747232,0.432867);  
rgb(652pt)=(0.614308,0.747129,0.431922);  
rgb(653pt)=(0.616708,0.747027,0.430978);  
rgb(654pt)=(0.619107,0.746925,0.430033);  
rgb(655pt)=(0.621487,0.746823,0.429089);  
rgb(656pt)=(0.623861,0.746721,0.428144);  
rgb(657pt)=(0.626235,0.746619,0.4272);    
rgb(658pt)=(0.628609,0.746517,0.426256);  
rgb(659pt)=(0.630962,0.746393,0.425311);  
rgb(660pt)=(0.63331,0.746266,0.424367);   
rgb(661pt)=(0.635658,0.746138,0.423422);  
rgb(662pt)=(0.638007,0.746011,0.422478);  
rgb(663pt)=(0.640332,0.745906,0.421557);  
rgb(664pt)=(0.642654,0.745804,0.420638);  
rgb(665pt)=(0.644977,0.745702,0.419719);  
rgb(666pt)=(0.6473,0.7456,0.4188);        
rgb(667pt)=(0.649623,0.745472,0.417881);  
rgb(668pt)=(0.651946,0.745345,0.416962);  
rgb(669pt)=(0.654268,0.745217,0.416043);  
rgb(670pt)=(0.656587,0.745089,0.415124);  
rgb(671pt)=(0.658859,0.744962,0.414205);  
rgb(672pt)=(0.661131,0.744834,0.413286);  
rgb(673pt)=(0.663402,0.744707,0.412368);  
rgb(674pt)=(0.665674,0.744579,0.411453);  
rgb(675pt)=(0.667946,0.744451,0.410559);  
rgb(676pt)=(0.670218,0.744324,0.409666);  
rgb(677pt)=(0.672489,0.744196,0.408773);  
rgb(678pt)=(0.674755,0.744062,0.407879);  
rgb(679pt)=(0.677001,0.743909,0.406986);  
rgb(680pt)=(0.679247,0.743756,0.406092);  
rgb(681pt)=(0.681494,0.743603,0.405199);  
rgb(682pt)=(0.68374,0.743458,0.404306);   
rgb(683pt)=(0.685986,0.74333,0.403412);   
rgb(684pt)=(0.688232,0.743203,0.402519);  
rgb(685pt)=(0.690479,0.743075,0.401626);  
rgb(686pt)=(0.692704,0.742937,0.400732);  
rgb(687pt)=(0.694899,0.742784,0.399839);  
rgb(688pt)=(0.697094,0.742631,0.398945);  
rgb(689pt)=(0.699289,0.742477,0.398052);  
rgb(690pt)=(0.701497,0.742324,0.397171);  
rgb(691pt)=(0.703718,0.742171,0.396303);  
rgb(692pt)=(0.705939,0.742018,0.395435);  
rgb(693pt)=(0.708159,0.741865,0.394568);  
rgb(694pt)=(0.710351,0.741712,0.3937);    
rgb(695pt)=(0.71252,0.741559,0.392832);   
rgb(696pt)=(0.71469,0.741405,0.391964);   
rgb(697pt)=(0.71686,0.741252,0.391096);   
rgb(698pt)=(0.719029,0.741082,0.390228);  
rgb(699pt)=(0.721199,0.740904,0.38936);   
rgb(700pt)=(0.723369,0.740725,0.388492);  
rgb(701pt)=(0.725538,0.740546,0.387625);  
rgb(702pt)=(0.727708,0.740386,0.386757);  
rgb(703pt)=(0.729878,0.740233,0.385889);  
rgb(704pt)=(0.732047,0.74008,0.385021);   
rgb(705pt)=(0.734217,0.739927,0.384153);  
rgb(706pt)=(0.736366,0.739753,0.383285);  
rgb(707pt)=(0.73851,0.739574,0.382417);   
rgb(708pt)=(0.740654,0.739395,0.38155);   
rgb(709pt)=(0.742798,0.739217,0.380682);  
rgb(710pt)=(0.744919,0.739038,0.379837);  
rgb(711pt)=(0.747038,0.738859,0.378995);  
rgb(712pt)=(0.749156,0.738681,0.378152);  
rgb(713pt)=(0.751275,0.738502,0.37731);   
rgb(714pt)=(0.753394,0.738323,0.376442);  
rgb(715pt)=(0.755512,0.738145,0.375574);  
rgb(716pt)=(0.757631,0.737966,0.374707);  
rgb(717pt)=(0.75975,0.737789,0.373841);   
rgb(718pt)=(0.761868,0.737636,0.372998);  
rgb(719pt)=(0.763987,0.737483,0.372156);  
rgb(720pt)=(0.766105,0.73733,0.371314);   
rgb(721pt)=(0.76822,0.737169,0.370471);   
rgb(722pt)=(0.770313,0.736965,0.369629);  
rgb(723pt)=(0.772406,0.73676,0.368786);   
rgb(724pt)=(0.774499,0.736556,0.367944);  
rgb(725pt)=(0.776592,0.736358,0.367096);  
rgb(726pt)=(0.778686,0.736179,0.366228);  
rgb(727pt)=(0.780779,0.736001,0.36536);   
rgb(728pt)=(0.782872,0.735822,0.364492);  
rgb(729pt)=(0.784957,0.735643,0.363632);  
rgb(730pt)=(0.787024,0.735465,0.36279);   
rgb(731pt)=(0.789092,0.735286,0.361948);  
rgb(732pt)=(0.791159,0.735107,0.361105);  
rgb(733pt)=(0.793227,0.734929,0.360263);  
rgb(734pt)=(0.795295,0.73475,0.359421);   
rgb(735pt)=(0.797362,0.734571,0.358578);  
rgb(736pt)=(0.79943,0.734392,0.357736);   
rgb(737pt)=(0.801485,0.734214,0.356881);  
rgb(738pt)=(0.803527,0.734035,0.356014);  
rgb(739pt)=(0.805569,0.733856,0.355146);  
rgb(740pt)=(0.807611,0.733678,0.354278);  
rgb(741pt)=(0.809668,0.733499,0.353424);  
rgb(742pt)=(0.811735,0.73332,0.352582);   
rgb(743pt)=(0.813803,0.733142,0.35174);   
rgb(744pt)=(0.81587,0.732963,0.350897);   
rgb(745pt)=(0.817921,0.732784,0.350038);  
rgb(746pt)=(0.819963,0.732606,0.349171);  
rgb(747pt)=(0.822005,0.732427,0.348303);  
rgb(748pt)=(0.824047,0.732248,0.347435);  
rgb(749pt)=(0.826071,0.73207,0.346567);   
rgb(750pt)=(0.828087,0.731891,0.345699);  
rgb(751pt)=(0.830104,0.731712,0.344831);  
rgb(752pt)=(0.83212,0.731534,0.343963);   
rgb(753pt)=(0.834158,0.731355,0.343095);  
rgb(754pt)=(0.8362,0.731176,0.342228);    
rgb(755pt)=(0.838242,0.730998,0.34136);   
rgb(756pt)=(0.840284,0.730819,0.340492);  
rgb(757pt)=(0.842303,0.73064,0.339624);   
rgb(758pt)=(0.84432,0.730462,0.338756);   
rgb(759pt)=(0.846336,0.730283,0.337888);  
rgb(760pt)=(0.848353,0.730104,0.33702);   
rgb(761pt)=(0.850369,0.729926,0.336153);  
rgb(762pt)=(0.852386,0.729747,0.335285);  
rgb(763pt)=(0.854402,0.729568,0.334417);  
rgb(764pt)=(0.856419,0.729391,0.333546);  
rgb(765pt)=(0.858435,0.729238,0.332627);  
rgb(766pt)=(0.860452,0.729085,0.331708);  
rgb(767pt)=(0.862468,0.728932,0.330789);  
rgb(768pt)=(0.864481,0.728778,0.329874);  
rgb(769pt)=(0.866472,0.728625,0.32898);   
rgb(770pt)=(0.868463,0.728472,0.328087);  
rgb(771pt)=(0.870454,0.728319,0.327194);  
rgb(772pt)=(0.872445,0.728166,0.326295);  
rgb(773pt)=(0.874436,0.728013,0.325376);  
rgb(774pt)=(0.876427,0.727859,0.324457);  
rgb(775pt)=(0.878418,0.727706,0.323538);  
rgb(776pt)=(0.880417,0.727561,0.322619);  
rgb(777pt)=(0.882433,0.727433,0.3217);    
rgb(778pt)=(0.88445,0.727306,0.320781);   
rgb(779pt)=(0.886466,0.727178,0.319862);  
rgb(780pt)=(0.888463,0.72705,0.318933);   
rgb(781pt)=(0.890429,0.726923,0.317989);  
rgb(782pt)=(0.892394,0.726795,0.317044);  
rgb(783pt)=(0.894359,0.726668,0.3161);    
rgb(784pt)=(0.896337,0.726552,0.315132);  
rgb(785pt)=(0.898328,0.72645,0.314136);   
rgb(786pt)=(0.900319,0.726348,0.313141);  
rgb(787pt)=(0.90231,0.726246,0.312145);   
rgb(788pt)=(0.904301,0.726158,0.31115);   
rgb(789pt)=(0.906292,0.726081,0.310154);  
rgb(790pt)=(0.908283,0.726005,0.309159);  
rgb(791pt)=(0.910274,0.725928,0.308163);  
rgb(792pt)=(0.912249,0.725851,0.307151);  
rgb(793pt)=(0.914214,0.725775,0.30613);   
rgb(794pt)=(0.91618,0.725698,0.305109);   
rgb(795pt)=(0.918145,0.725622,0.304088);  
rgb(796pt)=(0.920111,0.7256,0.303031);    
rgb(797pt)=(0.922076,0.7256,0.301959);    
rgb(798pt)=(0.924041,0.7256,0.300886);    
rgb(799pt)=(0.926007,0.7256,0.299814);    
rgb(800pt)=(0.927972,0.7256,0.298722);    
rgb(801pt)=(0.929938,0.7256,0.297624);    
rgb(802pt)=(0.931903,0.7256,0.296527);    
rgb(803pt)=(0.933869,0.7256,0.295429);    
rgb(804pt)=(0.935812,0.725668,0.294264);  
rgb(805pt)=(0.937752,0.725744,0.29309);   
rgb(806pt)=(0.939692,0.725821,0.291916);  
rgb(807pt)=(0.941632,0.725897,0.290741);  
rgb(808pt)=(0.943571,0.726023,0.289518);  
rgb(809pt)=(0.945511,0.726151,0.288293);  
rgb(810pt)=(0.947451,0.726278,0.287068);  
rgb(811pt)=(0.949389,0.726411,0.285839);  
rgb(812pt)=(0.951278,0.726641,0.284537);  
rgb(813pt)=(0.953167,0.72687,0.283235);   
rgb(814pt)=(0.955056,0.7271,0.281933);    
rgb(815pt)=(0.956938,0.72734,0.280622);   
rgb(816pt)=(0.958776,0.727646,0.279243);  
rgb(817pt)=(0.960614,0.727952,0.277865);  
rgb(818pt)=(0.962451,0.728259,0.276486);  
rgb(819pt)=(0.964273,0.728597,0.275086);  
rgb(820pt)=(0.966034,0.729057,0.273606);  
rgb(821pt)=(0.967795,0.729516,0.272126);  
rgb(822pt)=(0.969557,0.729976,0.270645);  
rgb(823pt)=(0.971288,0.730473,0.269135);  
rgb(824pt)=(0.972947,0.73106,0.267552);   
rgb(825pt)=(0.974606,0.731647,0.265969);  
rgb(826pt)=(0.976265,0.732234,0.264387);  
rgb(827pt)=(0.977857,0.732879,0.262785);  
rgb(828pt)=(0.979338,0.733619,0.261151);  
rgb(829pt)=(0.980818,0.734359,0.259518);  
rgb(830pt)=(0.982299,0.735099,0.257884);  
rgb(831pt)=(0.983697,0.73591,0.256227);   
rgb(832pt)=(0.984999,0.736803,0.254542);  
rgb(833pt)=(0.986301,0.737697,0.252858);  
rgb(834pt)=(0.987603,0.73859,0.251173);   
rgb(835pt)=(0.988753,0.739566,0.249474);  
rgb(836pt)=(0.989774,0.740613,0.247764);  
rgb(837pt)=(0.990795,0.741659,0.246054);  
rgb(838pt)=(0.991816,0.742706,0.244344);  
rgb(839pt)=(0.992677,0.743816,0.242681);  
rgb(840pt)=(0.993443,0.744965,0.241048);  
rgb(841pt)=(0.994209,0.746114,0.239414);  
rgb(842pt)=(0.994975,0.747262,0.23778);   
rgb(843pt)=(0.995578,0.748465,0.236165);  
rgb(844pt)=(0.996114,0.74969,0.234557);   
rgb(845pt)=(0.99665,0.750915,0.232949);   
rgb(846pt)=(0.997186,0.752141,0.231341);  
rgb(847pt)=(0.997562,0.753386,0.229813);  
rgb(848pt)=(0.997893,0.754637,0.228307);  
rgb(849pt)=(0.998225,0.755887,0.226801);  
rgb(850pt)=(0.998557,0.757138,0.225295);  
rgb(851pt)=(0.998711,0.758433,0.223856);  
rgb(852pt)=(0.998839,0.759735,0.222426);  
rgb(853pt)=(0.998966,0.761037,0.220997);  
rgb(854pt)=(0.999094,0.762339,0.219567);  
rgb(855pt)=(0.999076,0.763641,0.218186);  
rgb(856pt)=(0.99905,0.764942,0.216808);   
rgb(857pt)=(0.999025,0.766244,0.21543);   
rgb(858pt)=(0.998995,0.767546,0.214054);  
rgb(859pt)=(0.998868,0.768848,0.212752);  
rgb(860pt)=(0.99874,0.77015,0.21145);     
rgb(861pt)=(0.998613,0.771451,0.210149);  
rgb(862pt)=(0.998473,0.772756,0.208856);  
rgb(863pt)=(0.998243,0.774083,0.207631);  
rgb(864pt)=(0.998014,0.775411,0.206405);  
rgb(865pt)=(0.997784,0.776738,0.20518);   
rgb(866pt)=(0.997539,0.77806,0.20396);    
rgb(867pt)=(0.997232,0.779362,0.20276);   
rgb(868pt)=(0.996926,0.780664,0.201561);  
rgb(869pt)=(0.99662,0.781966,0.200361);   
rgb(870pt)=(0.996299,0.783268,0.199168);  
rgb(871pt)=(0.995942,0.784569,0.197994);  
rgb(872pt)=(0.995584,0.785871,0.19682);   
rgb(873pt)=(0.995227,0.787173,0.195646);  
rgb(874pt)=(0.994842,0.788475,0.19449);   
rgb(875pt)=(0.994408,0.789777,0.193367);  
rgb(876pt)=(0.993974,0.791078,0.192244);  
rgb(877pt)=(0.99354,0.79238,0.191121);    
rgb(878pt)=(0.993083,0.793671,0.190021);  
rgb(879pt)=(0.992598,0.794947,0.188949);  
rgb(880pt)=(0.992113,0.796223,0.187877);  
rgb(881pt)=(0.991628,0.797499,0.186805);  
rgb(882pt)=(0.99113,0.798789,0.185732);   
rgb(883pt)=(0.990619,0.800091,0.18466);   
rgb(884pt)=(0.990109,0.801393,0.183588);  
rgb(885pt)=(0.989598,0.802695,0.182516);  
rgb(886pt)=(0.989072,0.803996,0.18146);   
rgb(887pt)=(0.988536,0.805298,0.180413);  
rgb(888pt)=(0.988,0.8066,0.179367);       
rgb(889pt)=(0.987464,0.807902,0.17832);   
rgb(890pt)=(0.98691,0.809186,0.177291);   
rgb(891pt)=(0.986349,0.810462,0.17627);   
rgb(892pt)=(0.985787,0.811738,0.175249);  
rgb(893pt)=(0.985226,0.813015,0.174228);  
rgb(894pt)=(0.984644,0.814311,0.173207);  
rgb(895pt)=(0.984057,0.815613,0.172186);  
rgb(896pt)=(0.98347,0.816914,0.171165);   
rgb(897pt)=(0.982883,0.818216,0.170144);  
rgb(898pt)=(0.982296,0.819518,0.169145);  
rgb(899pt)=(0.981709,0.82082,0.16815);    
rgb(900pt)=(0.981122,0.822122,0.167154);  
rgb(901pt)=(0.980535,0.823423,0.166159);  
rgb(902pt)=(0.979947,0.824725,0.165163);  
rgb(903pt)=(0.97936,0.826027,0.164168);   
rgb(904pt)=(0.978773,0.827329,0.163172);  
rgb(905pt)=(0.978186,0.828631,0.162177);  
rgb(906pt)=(0.977599,0.829932,0.161181);  
rgb(907pt)=(0.977012,0.831234,0.160186);  
rgb(908pt)=(0.976425,0.832536,0.15919);   
rgb(909pt)=(0.975838,0.833841,0.158195);  
rgb(910pt)=(0.975251,0.835168,0.157199);  
rgb(911pt)=(0.974664,0.836495,0.156204);  
rgb(912pt)=(0.974077,0.837823,0.155208);  
rgb(913pt)=(0.973489,0.83915,0.154213);   
rgb(914pt)=(0.972902,0.840477,0.153217);  
rgb(915pt)=(0.972315,0.841805,0.152222);  
rgb(916pt)=(0.971728,0.843132,0.151226);  
rgb(917pt)=(0.971155,0.844466,0.150224);  
rgb(918pt)=(0.970619,0.845819,0.149203);  
rgb(919pt)=(0.970083,0.847172,0.148182);  
rgb(920pt)=(0.969547,0.848525,0.147161);  
rgb(921pt)=(0.96902,0.849886,0.14614);    
rgb(922pt)=(0.968509,0.851265,0.145119);  
rgb(923pt)=(0.967999,0.852643,0.144098);  
rgb(924pt)=(0.967488,0.854022,0.143077);  
rgb(925pt)=(0.967,0.855411,0.142056);     
rgb(926pt)=(0.966541,0.856815,0.141035);  
rgb(927pt)=(0.966081,0.858219,0.140014);  
rgb(928pt)=(0.965622,0.859623,0.138992);  
rgb(929pt)=(0.965189,0.86104,0.137945);   
rgb(930pt)=(0.96478,0.862469,0.136873);   
rgb(931pt)=(0.964372,0.863899,0.135801);  
rgb(932pt)=(0.963963,0.865328,0.134729);  
rgb(933pt)=(0.96357,0.866773,0.133657);   
rgb(934pt)=(0.963187,0.868228,0.132585);  
rgb(935pt)=(0.962805,0.869683,0.131513);  
rgb(936pt)=(0.962422,0.871138,0.130441);  
rgb(937pt)=(0.962091,0.87261,0.129351);   
rgb(938pt)=(0.961785,0.874091,0.128253);  
rgb(939pt)=(0.961478,0.875571,0.127156);  
rgb(940pt)=(0.961172,0.877052,0.126058);  
rgb(941pt)=(0.960885,0.878571,0.124961);  
rgb(942pt)=(0.960605,0.880103,0.123863);  
rgb(943pt)=(0.960324,0.881634,0.122765);  
rgb(944pt)=(0.960043,0.883166,0.121668);  
rgb(945pt)=(0.959849,0.884719,0.120549);  
rgb(946pt)=(0.95967,0.886276,0.119426);   
rgb(947pt)=(0.959491,0.887833,0.118302);  
rgb(948pt)=(0.959313,0.88939,0.117179);   
rgb(949pt)=(0.959181,0.890995,0.116032);  
rgb(950pt)=(0.959054,0.892603,0.114884);  
rgb(951pt)=(0.958926,0.894211,0.113735);  
rgb(952pt)=(0.958799,0.895819,0.112587);  
rgb(953pt)=(0.958748,0.897453,0.111464);  
rgb(954pt)=(0.958697,0.899086,0.110341);  
rgb(955pt)=(0.958646,0.90072,0.109217);   
rgb(956pt)=(0.958602,0.902359,0.108089);  
rgb(957pt)=(0.958628,0.904043,0.106915);  
rgb(958pt)=(0.958653,0.905728,0.105741);  
rgb(959pt)=(0.958679,0.907413,0.104567);  
rgb(960pt)=(0.958718,0.909102,0.103393);  
rgb(961pt)=(0.95882,0.910812,0.102219);   
rgb(962pt)=(0.958922,0.912522,0.101044);  
rgb(963pt)=(0.959024,0.914232,0.0998703); 
rgb(964pt)=(0.959153,0.915962,0.0986961); 
rgb(965pt)=(0.959357,0.917749,0.0975219); 
rgb(966pt)=(0.959561,0.919536,0.0963477); 
rgb(967pt)=(0.959765,0.921323,0.0951736); 
rgb(968pt)=(0.959996,0.923118,0.0939907); 
rgb(969pt)=(0.960277,0.924931,0.092791);  
rgb(970pt)=(0.960557,0.926743,0.0915913); 
rgb(971pt)=(0.960838,0.928555,0.0903916); 
rgb(972pt)=(0.961151,0.930378,0.0891919); 
rgb(973pt)=(0.961509,0.932216,0.0879922); 
rgb(974pt)=(0.961866,0.934054,0.0867925); 
rgb(975pt)=(0.962223,0.935892,0.0855928); 
rgb(976pt)=(0.96262,0.937768,0.0843802);  
rgb(977pt)=(0.963053,0.939683,0.083155);  
rgb(978pt)=(0.963487,0.941597,0.0819297); 
rgb(979pt)=(0.963921,0.943512,0.0807045); 
rgb(980pt)=(0.964415,0.945426,0.0794643); 
rgb(981pt)=(0.964951,0.947341,0.0782135); 
rgb(982pt)=(0.965487,0.949255,0.0769628); 
rgb(983pt)=(0.966023,0.951169,0.075712);  
rgb(984pt)=(0.966594,0.953118,0.0744441); 
rgb(985pt)=(0.967181,0.955083,0.0731679); 
rgb(986pt)=(0.967768,0.957049,0.0718916); 
rgb(987pt)=(0.968355,0.959014,0.0706153); 
rgb(988pt)=(0.96898,0.960999,0.0693006);  
rgb(989pt)=(0.969619,0.96299,0.0679733);  
rgb(990pt)=(0.970257,0.964981,0.0666459); 
rgb(991pt)=(0.970895,0.966972,0.0653186); 
rgb(992pt)=(0.971554,0.968984,0.0639486); 
rgb(993pt)=(0.972218,0.971001,0.0625703); 
rgb(994pt)=(0.972882,0.973017,0.0611919); 
rgb(995pt)=(0.973545,0.975034,0.0598135); 
rgb(996pt)=(0.974232,0.97705,0.058318);   
rgb(997pt)=(0.974922,0.979067,0.056812);  
rgb(998pt)=(0.975611,0.981083,0.055306);  
rgb(999pt)=(0.9763,0.9831,0.0538)
}      
}   

%% file: usermacros.tex
\usepgfplotslibrary{groupplots}

\newlength{\figureheight}
\newlength{\figurewidth}
\usetikzlibrary{arrows,shapes,positioning}
\usetikzlibrary{decorations.markings}
\tikzstyle arrowstyle=[scale=1]
\tikzstyle directed=[postaction={decorate,decoration={markings,
		mark=at position .65 with {\arrow[arrowstyle]{stealth}}}}]
\tikzstyle reverse directed=[postaction={decorate,decoration={markings,
		mark=at position .65 with {\arrowreversed[arrowstyle]{stealth};}}}]

\pgfplotscreateplotcyclelist{color parula}{% 
	royalblue!70,every mark/.append style={solid,line width = 0pt,fill=royalblue!60!black},mark=ball\\%
	goldenrod!50!yelloworange,every mark/.append style={solid,fill=goldenrod!30!black},mark=*\\%
	green,every mark/.append style={black,fill=yellowgreen},mark=diamond*\\%
	bluegreen,every mark/.append style={fill=black},mark=triangle*\\%	
	royalblue!70,densely dashed,every mark/.append style={solid,line width = 0pt,fill=royalblue!60!black},mark=ball\\%
	goldenrod!50!yelloworange,densely dashed,every mark/.append style={solid,black,fill=goldenrod},mark=diamond*\\%
	yellowgreen,densely dashed,every mark/.append style={solid,fill=yellowgreen!60!royalblue},mark=square*, mark size= 1pt\\%
	bluegreen,densely dashed,mark size = 1.5pt,every mark/.append style={solid,fill=white},mark=otimes*\\%,densely dashed	
	royalblue,densely dashed,mark size = 1.5pt,every mark/.append style={solid,fill=royalblue!30!black},mark=pentagon*\\%,densely dashed
	goldenrod!40!yelloworange,every mark/.append style={solid,fill=goldenrod!30!black},mark=|\\%		
	}

%% file: Sections/introduction.tex
\section{Introduction}
Population balance equations are widely used in engineering applications, 
including aerosol physics, high shear granulation,  pharmaceutical industries,  polymerization and emulsion processes,   evaporation and condensation processes in bubble column reactors, bioreactors, turbulent flame reactors and many others, see \cite{RANDOLPH1971,Ramkrishna2000,Rosner2003,Cameron2005,Dhanasekharan2005,Gokhale2009,Gokhale2009} and references therein.
These polydisperse processes are characterized by two phases: one of them is the continuous
and the second is a dispersed phase consisting of particles. The particles can take different forms like
crystals, drops or bubbles with several possible properties such as volume, chemical composition, porosity and enthalpy. In this work only the volume is considered.\\
The dynamic evolution of the particle number distribution which is described by the population balance equation (PBE) of the dispersed phase depends not only on the particle-particle interactions, but also on the continuous phase, due to interaction of these particles with the continuous flow field in which they are dispersed.
These interactions usually result in the common mechanisms aggregation, breakage, condensation, growth and nucleation. 
In our work we concentrate on binary aggregation, namely the Smoluchowski aggregation \cite{Qamar2007} and multiple breakages, since binary breakage is not sufficient for some of these applications.
Binary aggregation is the process of merging two particles to a larger particle, whereas in a breakage process, a particle breaks into several smaller fragments. 
Often included in a typical PBE are spatial transport terms, i.e. advection and diffusion terms. We use the simplest case and assume that the mean particle velocity is the same as those of the fluid.
The resulting population balance equations range from integro-differential to partial integro-differential equations of hyperbolic or parabolic type in phase space, in addition to the differential terms for advection and diffusion in the physical space.\\
The complex structure of the PBE allows for an analytical solution only for some simple breakage and aggregation kernels, see \cite{Dubouskii1992,Ziff1991} and references therein. Thus, numerical solution methods have been extensively studied.
In literature, several different classes of methods are discussed like  Monte-Carlo methods, see, for example  \cite{Lee2000} or finite element techniques, see for example \cite{Mahoney2002} and references therein.
Other numerical techniques are the method of successive approximations \cite{Ramkrishna2000}, the
methods of classes \cite{Ramkrishna2000,Kumar2006,Vanni2000}, and the method of weighted residuals \cite{Ramkrishna2000,Wulkow2001}.
Another approach is to apply finite volume schemes (FVS), which are frequently used for solving conservation laws. In \cite{Filbet2004a,Filbet2004}, this approach is used for solving aggregation PBEs using a conservative formulation of the equation. Later, the scheme has been applied to solve the combined breakage and aggregation PBE in  \cite{Kumar2009}. In  \cite{Kumar2012}, sectional methods
have been developed, among them the cell average technique for both uniform and non-uniform grids.\\
Another class of methods for solving the complex PBE is the method of moments. Instead of solving the PBE in full phase space, moment models reduce the high dimensionality of the original equation by using a finite set of moment equations, which are obtained by taking volume averages with respect to some basis of the volume space. The resulting moment equations require additional information about
the unknown particle number distribution that must be approximated via a closure. By the use of such an ansatz function the distribution function can be reconstructed from the available moments. 
One big issue regarding the method of moments is realizability. Roughly speaking, realizable moments are those moments associated with non-negative distribution functions. The set of realizable moments forms a convex cone in the set of all moment vectors.
Several types of  closures are  available. One that is often used for population balance equations is the quadrature method of moments ($QMOM$), first introduced  in the PBE context in \cite{McGraw1997}. It uses a non-linear Gaussian-like quadrature rule to compute weights and abscissas from the available moments, assuming  that the underlying distribution function is a weighted sum of Dirac $\delta-$functions. The computation of the quadrature rule involves solving a non-linear system. There also exist several extensions of $QMOM$, such as $SQMOM$ \cite{Attarakih2009,DTAKB08,AJDBTSK09} and $EQMOM$ \cite{Yuan2012}.\\
The standard moment method in linear radiative transport is the polynomial closure $P_N$ \cite{PNintro}, where $N$ is the order of the highest-order moment. For multi-variable cell population balance equations, $P_N$ closures of different types where discussed  in \cite{Mantzaris2001}. Using the polynomial closure with respect to an orthogonal basis ensures a non-expensive moment method, at the price of a potentially negative, and thus physically meaningless, number density.  
Another moment approach is the use of entropy-based closures $M_N$ which approximate the full distribution by an ansatz that solves a constrained, convex optimization problem. Using physically relevant entropies, the big advantage of $M_N$ is the preservation of many important properties like positivity, entropy dissipation, and hyperbolicity \cite{Levermore1996}.
The disadvantage is the computational cost, especially compared to the polynomial closure, since an optimization problem has to be solved numerically in every point of the space-time grid. Realizability plays a big role for $M_N$, since the solving the optimization problem near the realizability boundary is especially expensive or even impossible.\\
In the present paper we aim at a numerical comparison of the above described closure approaches. 
The outline of the paper is the following. First, in section \ref{PBEsec} we introduce the population balance equation used in this paper. Then we shortly state the finite volume scheme, based on the work in  \cite{Kumar2014}. After that in section \ref{MomentModels} we explain the method of moments and the concept of realizability. Also, contained in this section are the different types of moment closures for which we compute the PBE, namely the polynomial closure $P_N$, the maximum entropy closure $M_N$ and the quadrature method of moments $QMOM_N$. This is followed by a description of the numerical  methods, supplemented with further implementation details in section \ref{implementation}. The numerical results for the  homogeneous case neglecting diffusion and advection are presented in section \ref{results}. These results are followed by  two-dimensional results on a PBE coupled to fluid dynamic equations for a lid-driven cavity test.
In section \ref{outlook} we summarize our conclusions.

%% file: Sections/PBE.tex
\section{Population Balance Equation}
\label{PBEsec}
We consider  a one-dimensional population balance equation (PBE). It  describes the time evolution of the particle number distribution function $f(t,x,v)$ under the simultaneous effect of binary aggregation, multiple breakage, diffusion and advection processes \cite{Kumar2014}.  The PBE reads
 \begin{equation}
 \begin{aligned}
  \label{PBE}
  \partial_t f(t,x&,v)+\nabla_x\cdot(v_d(t,x)f(t,x,v))-\nabla_x\cdot(D(t,x)\nabla_x f(t,x,v) )\\ =&\underbrace{\frac{1}{2}\int_{v_{min}}^{v} \omega(t,x,v-v',v')f(t,x,v-v')f(t,x,v')dv'}_\text{$=:A^+[f](t,x,v)$}-\underbrace{f(t,x,v)\int_{v_{min}}^{v_{max}}\omega(t,x,v,v')f(t,x,v')dv'}_\text{$=:A^-[f](t,x,v)$}\\
&+\underbrace{\int_v^{v_{max}}\Gamma (t,x,v')\beta (v,v')f(t,x,v')dv'}_\text{$=:B^+[f](t,x,v)$}-\underbrace{\Gamma(t,x,v)f(t,x,v)}_\text{$=:B^-[f](t,x,v)$}.
\end{aligned}
\end{equation}
The particle number density $f(t,x,v)\geq 0$ depends on the time $t\in\mathbb{R}^+$, the space $x\in\mathbb{R}^{n_{d}},~n_{d}\geq1$ and the volume $v\in[v_{min},v_{max}]\subseteq\mathbb{R}^+$ of the particles. The particle volume is limited through the minimal volume $v_{min}$ and the maximal volume $v_{max}$, which are physical properties determined by experiments. In theoretical works, $v_{min}=0$ and $v_{max}=\infty$ is widely used. The left-hand side of the equation is devoted to the space and time evolution of the density: it includes an advection term with the droplet velocity $v_d(t,x)$ and a diffusion term including a diffusion rate $D(t,x)$.
On the right hand side $A^+$ and $A^-$ characterize the binary aggregation, where the birth term $A^+$ determines the creation of particles of volume $v$ and the death term $A^-$ describes disappearance of particles of volume $v$ in the population balance.  $B^+$ and $B^-$ are the breakage terms \cite{Kumar2014}. The  aggregation and breakage processes are modeled via  the aggregation kernel $\omega(t,x,v,v')$, the breakage frequency $\Gamma(t,x,v)$ and the daughter droplet distribution function $\beta(v,v').$
The breakage frequency $\Gamma(t,x,v)\geq0$ represents the fractional number of breakage events per unit time of droplets of size $v$ while the daughter droplet distribution function $\beta(v,v')\geq0$ is the probability function for the creation of particles of size $v$ from particles of size $v'$ \cite{Kumar2006a}.

We assume the following properties for the daughter droplet distribution function:
\begin{subequations}
\label{B}
\begin{align}
 \label{B1}
  \beta(v,v')&\equiv0,~~~\forall v'<v,\\
 \label{B2}
 \int_{v_{min}}^{v'}\beta(v,v')dv&=N(v'),~~\forall v'\in[v_{min},v_{max}],\\
 \label{B3}
 \int_{v_{min}}^{v'}v\beta(v,v')dv&=v',~~\forall v'\in[v_{min},v_{max}].
 \end{align}
\end{subequations}
The first property \eqref{B1} comes from the fact that no particle can split into larger particles. The function $N(v')\in\mathbb{R}^+$ in \eqref{B2} represents the number of fragments obtained from the splitting of a particle of size $v'$ \cite{Kumar2006a}. Property \eqref{B3} takes into account that over time the total mass should be conserved by breakage events, while the total number of particles increases.

We furthermore set $\beta(v,v')\equiv 0,~\Gamma(t,x,v)\equiv0$ for all $v,v'\notin[v_{min},v_{max}]$.

Binary aggregation appears if two particles of size $v_1$ and $v_2$ collide and merge to a particle of size $v_1+v_2$. The aggregation kernel $\omega$ provides the probability that a collision of two particles results in a successful merge. We again assume the following conditions:
 \begin{subequations}
 \label{A}
 \begin{align}
 \label{A0}
  \omega(t,x,v,v')&\geq0,~~~\forall v,v'\in[v_{min},v_{max}],\\
 \label{A1}
   \omega(t,x,v,v')&=\omega(t,x,v',v),~~~\forall v,v'\in[v_{min},v_{max}],\\
   \label{A2}
    \omega(t,x,v,v')&\equiv0,~~~\forall v,v' \in[v_{min},v_{max}]\text{ with } v+v'>v_{max},\\
 \label{A3}
    \omega(t,x,v,v')&\equiv 0,~~~\forall v,v'\notin[v_{min},v_{max}].
         \end{align}
 \end{subequations}
The symmetry condition \eqref{A1} comes from the fact that the two merging events $(v_1,v_2)\to v_1+v_2$ and $(v_2,v_1)\to v_1+v_2$ are equivalent (symmetry of binary aggregation events). 

Formulas for our choice of kernels can be found in section \ref{results}.

%% file: Sections/momentmodels.tex
\section{Moment Models and Realizability}
\label{MomentModels}
 The method of moments is commonly used to derive reduced models for kinetic transport equations \cite{Schneider2014,Yuan2011,Garrett2013,Hauck2011,Vikas2013a}.
 To solve such a reduced model is often less expensive, regarding computational time, than to solve the complete model. 
 Equation \eqref{PBE} depends on  time, space and the volume. Here,  the method of moments consists in a (nonlinear) projection of the solution $f$ with respect to a polynomial basis in the volume variable $v$. 
 
We start with some definitions that closely follow \cite{Levermore1996, Alldredge2012, Alldredge2014}.

We denote by $\gamma=\left(\gamma_0,\gamma_1,...,\gamma_N\right)^T:\mathbb{R}^+\times\mathbb{R}^{n_d}\to\mathbb{R}^{N+1}$ the vector containing entries of the first $N+1$ moments of $f$ with respect to some polynomial basis $m=m(v)=\left(m_0(v),...,m_N(v)\right)^T$ of $\mathbb{P}_N([v_{min},v_{max}])$\footnote{Other bases are also possible, see, e.g., \cite{Schneider2014}.}, that means
\begin{equation}
\begin{aligned}
\label{Momente}
\gamma(t,x):=\langle mf(t,x,\cdot)\rangle := \int_{v_{min}}^{v_{max}}m(v)f(t,x,v)~dv,
\end{aligned}
\end{equation}
where the integration is applied component-wise.
Then the system of $N+1$ moment equations for the population balance equation is derived by multiplying \eqref{PBE} by the basis $m$ and integrating over the volume domain. This results in
\begin{equation}
\begin{aligned}
\label{ME}
\partial_t \gamma(t,&x)+\nabla_x\cdot(v_d(t,x)\gamma(t,x))-\nabla_x\cdot(D(t,x)\nabla_x \gamma(t,x) )\\ =&\langle m(v)A^+[f](t,x,v) \rangle-\langle m(v)A^-[f](t,x,v) \rangle+\langle m(v)B^+[f](t,x,v) \rangle-\langle m(v)B^-[f](t,x,v) \rangle.  
\end{aligned}
\end{equation}
The moments $\gamma_i$ depend on  the distribution function that appears on the right hand side of \eqref{ME} in the aggregation and breakage terms.
Here,  an  ansatz for the distribution function, which depends on the known moments, is needed to close  the equations. 

\subsection{Realizability}
When dealing with closures, the question arises which moments can actually be realized by a physical (non-negative) distribution function, compare  \cite{Curto1991, REA2}. 
In other words, for given moments $\gamma_i$ one has  to find  a non-negative distribution function $f(v)\geq 0$ such that
$$\langle m_i(v)f(v)\rangle=\gamma_i,~~~0\leq i\leq N.$$

The realizability domain $$\textbf{\textit{R}}_N = \{\gamma\in\mathbb{R}^{N+1}~|~\exists f\geq 0,~ \langle mf\rangle = \gamma\}$$ is defined as the set of vectors $(\gamma_0,...,\gamma_N)^T\in\mathbb{R}^{N+1}$ such that this problem has a solution.
We refer to \cite{Curto1991, REA2} for further reading and explicit characterizations of $\textbf{\textit{R}}_N$ in terms of the moments $\gamma$.
Note that the realizable set changes when $\langle\cdot\rangle$ is replaced by a numerical quadrature, see section \ref{Mereal} and \cite{Alldredge2014,Alldredge2015}

\subsection{Moment Closures}
Assuming a specific form of the distribution function $f$ in  the integrals appearing in \eqref{ME}, which depends on the available moments, one can close the equations.
There exist many types of closures in literature, each with their own advantages and disadvantages. In this section we will explain the polynomial closure $P_N$, the maximum entropy closure $M_N$ as well as the quadrature method of moments $QMOM$. %The $P_N$ method does not guarantee positivity of the underlying distribution function. By definition the maximum entropy models and $QMOM$ use non-negative distribution functions to close the equations. However, in contrast to the maximum entropy model,  for a discretization of the $QMOM$ moment equation, one is not able to predict a time step size such that the moments remain realizable in every time step.
%On the other hand, in case of the maximum entropy models, the solution the  optimization problem in every time step  leads to  higher computational costs than solving the $P_N$ moment equations. 

\subsubsection{Polynomial Closure}
The $P_N$ equations can be most easily understood as a Galerkin semi-discretization in $v$. The idea is to expand the distribution function in terms of a truncated series expansion 
$$f_{P_N}(t,x,v)=\sum_{i=0}^Na_i(t,x)m_i(v).$$
Inserting $f_{P_N}$ into the definition of moments \eqref{Momente} leads to
 \begin{equation}
 \begin{aligned}
  \label{li}
\gamma_i(t,x)=\sum_{j=0}^Na_j(t,x)\langle m_i(v)m_j(v)\rangle,~~~i=0,...,N.
\end{aligned}
\end{equation}
The basis $m$ can be chosen to be any kind of polynomial basis, for example the monomials $v^i$. But to simplify the solution of the linear system of equations \eqref{li} for the expansion coefficients $a_i$, we use an orthogonal basis. Choosing the Legendre polynomials and shifting them to the interval $[v_{min},v_{max}]$ gives us the following recursion formula 
\begin{equation}
\begin{aligned}
\label{Legendrebasis}
m_0(v)=1,~m_1(v)=\frac{2}{v_{max}-v_{min}}&(v-v_{min})-1,\\
m_i(v)=\frac{2i-1}{i}\left(\frac{2}{v_{max}-v_{min}}(v-v_{min})-1\right)&m_{i-1}(v)-\frac{i-1}{i}m_{i-2}(v), \text{ for }i\geq2.
\end{aligned}
\end{equation}
Therefore, the moments $\gamma_i$ for $i=0,...,N$ can be expressed as
\begin{align*}
\gamma_i(t,x)&=\sum_{j=0}^Na_j(t,x)\langle m_im_j\rangle=\sum_{j=0}^Na_j(t,x)\tfrac{v_{max}-v_{min}}{2}\tfrac{2}{2i+1}\delta_{ij}=a_i(t,x)\tfrac{v_{max}-v_{min}}{2i+1},
\end{align*}
where $\delta_{ij}$ is the Kronecker delta. Solving for the coefficients $a_i$, the approximated distribution function for the $P_N$ closure looks like 
\begin{equation}
\begin{aligned}
\label{PN}
f_{P_N}(t,x,v)=\sum_{i=0}^N\gamma_i(t,x)\frac{2i+1}{v_{max}-v_{min}}m_i(v).
\end{aligned}
\end{equation}
The polynomial closure does not necessarily provide a non-negative distribution function, thereby the computed moments could get physically meaningless, like e.g. a negative particle number.
\subsubsection{Maximum Entropy}
 Another kind of closure for the  approximation of the distribution $f$ are  entropy-based methods. 
 Here, $f_{M_N}$ is the solution of a constrained, convex optimization problem \cite{Hauck2011}.
 Namely, for strictly convex $\eta:\mathbb{R}\to\mathbb{R}$ the optimization problem looks like
  \begin{equation}
     \begin{aligned}
  \label{max1}
  \min_{g}&\langle\eta(g)\rangle,\\
 \text{subject to }&\langle mg\rangle=\gamma.
     \end{aligned}
\end{equation}
 Transforming \eqref{max1} into its unconstrained, strictly convex dual problem we end at searching for the Lagrange multipliers $\hat{\alpha}(\gamma)\in\mathbb{R}^{N+1}$ that solve
 $$\min_{\alpha\in\mathbb{R}^{N+1}}\left(\langle\eta(g)+\alpha^Tmg\rangle-\alpha\gamma\right).$$
Using the Legendre dual $\eta_*:\mathbb{R}\to\mathbb{R}$ of $\eta$, this is equivalent to
 \begin{equation}
 \begin{aligned}
  \label{max}
   \hat{\alpha}(\gamma)=\underset{{\alpha\in\mathbb{R}^{N+1}}}{\text{argmin}}\left(\langle\eta_*(\alpha^Tm)\rangle-\alpha\gamma\right).
    \end{aligned}
\end{equation}
 So, if a solution of problem \eqref{max} exists, it takes the form \cite{Levermore1996}
  \begin{equation}
     \begin{aligned}
  \label{max2}
  f_{M_N}=G_{\hat{\alpha}}:=\eta_*'(\hat{\alpha}^Tm).
  \end{aligned}
  \end{equation}
 Choosing the Maxwell-Boltzmann entropy $\eta(g)=g\log(g)-g$, the Legendre dual and its derivative are of the form $\eta_*(y)=\eta_*'(y)=e^y$
 and therefore the solution of \eqref{max} can be written as
  $$G_{\alpha}=\exp(\alpha^Tm).$$ 
We define the dual objective function $h:\mathbb{R}^{N+1}\to\mathbb{R}$ (the function which gets minimized) by
  $$h(\alpha)=\langle G_{\alpha}\rangle-\alpha^T\gamma.$$
 Through the properties of the exponential function it is ensured that $f_{M_N}$ is always a positive function. Because of the strict convexity, problem \eqref{max1} does have a unique solution whenever $\gamma$ is realizable \cite{Hauck2011}.
%Moreover, realizability can be  guaranteed in every time step for a suitable time discretization of the moment system,
%see Proposition \ref{propm} below.

\subsubsection{QMOM}
In $QMOM$ the general idea is to use an n-atomic discrete measure as ansatz. This means that the approximated distribution function consists only of a linear combination of Dirac $\delta$-functions 
  \begin{equation}
     \begin{aligned}
  \label{QMOMf}
  f_{QMOM_N}=\sum_{j=1}^nw_j\delta(v-v_j),
  \end{aligned}
  \end{equation}
with $w_j>0$ and $v_j\in[v_{min},v_{max}]$ for all $j=1,...,n.$
Plugging $f_{QMOM_N}$ into \eqref{Momente} and choosing $m_i(v)=v^i$ leads to 
$$\gamma_i=\sum_{j=1}^nw_jv_j^i,~~~i=0,...,N.$$
The above nonlinear system can be solved using the Wheeler-algorithm \cite{C.Wheeler1974,Yuan2012,Curto1991,REA2} which diagonalizes a tridiagonal matrix to find the weights $w_i$ and abscissas $v_i$. This results in a robust and efficient algorithm for the inversion of the moment problem, which will only succeed if $\gamma$ is realizable. 
It can be shown that a moment vector on the realizability boundary can be uniquely represented by an atomic distribution function. Thus, $QMOM$ is able to exactly reproduce this behavior in such a  case \cite{Vikas2013a, Kershaw1976}.

%% file: Sections/modelling.tex
\section{Numerical Implementation}
\label{implementation}

\subsection{Finite Volume Scheme for the PBE}
\label{FVSchap}
We will shortly recall the finite volume scheme (FVS) for the one-dimensional PBE with binary aggregation and multiple breakages on uniform meshes, as introduced in \cite{Filbet2004,Bourgade2007}. Its solution will be used as a reference for our numerical comparison of the moment closure schemes. 

Extensions of this scheme to different types of uniform and non-uniform grids and analytical results can be found in \cite{Kumar2011,Kumar2014,Kumar2006}.

At first, equation \eqref{PBE} is rewritten in conservative form \cite{Kumar2014} for the mass density $vf(t,v)$,
$$\partial_t(vf(t,x,v))+\nabla_x\cdot(v_d(t,x)f(t,x,v))-\nabla_x\cdot(D(t,x)\nabla_x f(t,x,v) )+\partial_v\left(F^{agg}(t,x,v)+F^{brk}(t,x,v)\right)=0,$$
with the continuous aggregation flux
$$F^{agg}(t,x,v)=\int_{v_{min}}^v\int_{v-u}^{v_{max}}u\omega(t,x,u,w)f(t,x,u)f(t,x,w)dwdu,$$
and the continuous breakage flux
$$F^{brk}(t,x,v)=-\int_{v}^{v_{max}}\int_{v_{min}}^vu\beta(u,w)\Gamma(t,x,w)f(t,x,w)dudw.$$
The volume domain $[v_{min},v_{max}]$ is discretized into equidistant cells 
$$\Lambda_i:=[v_{i-\frac{1}{2}},v_{i+\frac{1}{2}}[,~~~i=1,...,n_v \text{ with }v_{i+\frac{1}{2}}=v_{min}+i\Delta v.$$
The grid size is denoted by $\Delta v=(v_{max}-v_{min})/n_v$ and the midpoints $v_i$ are computed by $v_i=\left(v_{i+\frac{1}{2}}+v_{i-\frac{1}{2}}\right)/2$.
Let $f_i(t)$ be the semi-discrete approximation of the cell average of the solution on cell $i$ \cite{Kumar2011}
$$f_i(t,x)=\frac{1}{\Delta v }\int_{v_{i-\frac{1}{2}}}^{v_{i+\frac{1}{2}}}f(t,x,v)dv.$$
The FVS is derived by integrating the conservation law over every cell $\Lambda_i,$
$$\frac{df_i(t,x)}{dt}=-\frac{1}{\Delta v v_i}\left(J_{i+\frac{1}{2}}^{agg}(t,x)-J_{i-\frac{1}{2}}^{agg}(t,x)+J_{i+\frac{1}{2}}^{brk}(t,x)-J_{i-\frac{1}{2}}^{brk}(t,x)\right).$$
The numerical fluxes $J_{i+\frac{1}{2}}$ are chosen as appropriate approximations of the continuous flux functions  \cite{Kumar2011,Kumar2012,Kumar2014,Kumar2006a,Filbet2004,Filbet2004}. 
The numerical breakage flux takes the form \cite{Kumar2014}
\begin{align*}
J_{i+\frac{1}{2}}^{brk}(t,x)=-\sum_{l=i+1}^{n_v}f_l(t)\Gamma(t,x,v_l)\Delta v^2\sum_{j=1}^{i}v_j\beta(v_j,v_l).
\end{align*}
Moreover, $J_{1/2}^{brk}(t,x)=J_{n_v+1/2}^{brk}(t,x)=0$ since at the boundary we have $F^{brk}(t,x,v_{min})=F^{brk}(t,x,v_{max})=0.$

The numerical aggregation flux is given by \cite{Filbet2004,Kumar2014}
\begin{align*}
J_{i+\frac{1}{2}}^{agg}(t,x)=&\sum_{l=1}^{i-\zeta}v_lf_l(t,x)\Delta vf_{i-l+1-\zeta}(t,x)\omega\left(t,x,v_{i-l+1-\zeta},v_l\right)\left(v_{min}+\left(\tfrac{1}{2}-\zeta\right)\Delta v\right)
\\&+\sum_{l=1}^{i-\zeta}v_lf_l(t,x)\Delta v\sum_{j=i-l+2-\zeta}^{n_v}f_j(t,x)\omega(t,x,v_j,v_l)\Delta v
\\&+\sum_{l=\max\left(i-\zeta+1,1\right)}^{i}v_lf_l(t,x)\Delta v\sum_{j=1}^{n_v}f_j(t,x)\Delta v\omega(t,x,v_j,v_l)
\end{align*}
with $\alpha_{i,l}=i-l-\zeta+2$. The index $\zeta$ is determined by 
$v_{i+\frac{1}{2}}-v_{l}=v_{i-l+1}-v_{min}\in\Lambda_{i-l+1-\zeta}$ 
and
$$\zeta=\left \lceil{\frac{v_{min}}{\Delta v}-\frac{1}{2}}\right \rceil.$$
Again, $J_{1/2}^{agg}(t,x)=J^{agg}_{n_v+1/2}(t,x)=0$, since at the boundaries we have $F^{agg}(t,x,v_{min})=F^{agg}(t,x,v_{max})=0$.
We remark that  restrictions of the size of  the time-step to ensure positivity of the distribution function have been investigated in  \cite{Filbet2004, Kumar2006a}.

\subsection{Numerical approximation of the moment system}

 We reconsider  the moment equations \eqref{ME}.
For a better numerical approximation we transform the second term on the right hand side of \eqref{ME} towards
\begin{align*}
\langle m(v)A^-[f](t,v) \rangle&\overset{\eqref{A3}}{=}\int_{v_{min}}^{v_{max}}m(v)f(t,v)\int_{v_{min}}^{v_{max}-v+v_{min}}\omega(t,v,v')f(t,v')dv'dv\\
&\hspace{0.22cm}=\int_{v_{min}}^{v_{max}}m(v)f(t,v)\int_{v_{min}}^{v_{max}}\omega(t,v,g_{A^-}(v,v'))f(t,g_{A^-}(v,v'))g_{A^-}'(v,v')dv'dv,
\end{align*}
with $g_{A^-}(v,x)=(v_{max}-v)(x-v_{min})/(v_{max}-v_{min})+v_{min}$. In general, one is not able to compute the above integrals  exactly. They have to be evaluated   numerically. As quadrature rule for the computation of the terms on the right hand side of \eqref{ME} we use a Gau{\ss} Lobatto formula, which integrates polynomials up to degree $2n_Q-3$ exactly for $n_Q$ quadrature points.\\
For the $QMOM_N$, whose closure is based on Dirac $\delta$ functions, we approximate the integrals by
\begin{align*}
&\langle m(v)A^+[f_{QMOM_N}]\rangle=\frac{1}{2}\sum_{i=1}^n\sum_{j=1}^nm(v_i+v_j)w_iw_j\omega(v_i,v_j),\\
&\langle m(v)A^-[f_{QMOM_N}]\rangle=\sum_{i=1}^n\sum_{j=1}^nm(v_j)w_iw_j\omega(v_j,v_i),\\
&\langle m(v)B^+[f_{QMOM_N}]\rangle=\sum_{i=1}^nw_i\Gamma(v_i)\int_{v_{min}}^{v_i}m(v)\beta(v,v_i)dv,\\
&\langle m(v)B^-[f_{QMOM_N}]\rangle=\sum_{i=1}^nm(v_i)w_i\Gamma(v_i).
\end{align*}
Using property \eqref{B3} of the $\beta-$kernel, we can easily see that mass conservation is exactly fulfilled in case of the quadrature moment method, when $m$ is the monomial basis:
$$\int_{v_{min}}^{v_{max}}v\left(B^+[f_{QMOM_N}](t,v)-B^-[f_{QMOM_N}](t,v)+A^+[f_{QMOM_N}](t,v)-A^-[f_{QMOM_N}](t,v)\right)dv=0.$$ 
 
\subsection{Numerical Realizability}
\label{Mereal}
Since for maximum-entropy models almost all of the integrals $\langle \cdot\rangle$ have to be evaluated numerically, we need to adapt the definition of the realizable set. 

For a function $g:[v_{min},v_{max}]\to\mathbb{R}$ the nodes $\left\{v_i\right\}_{i=1}^{n_Q}$ and weights $\left\{w_i\right\}_{i=1}^{n_Q}$ of a quadrature rule $\textbf{\textit{Q}}$ are chosen such that $\langle g \rangle$ is approximated by 
\begin{equation}
\begin{aligned}
\label{quad}
\langle g \rangle \approx \textbf{\textit{Q}}(g)=\sum_{i=1}^{n_Q}w_ig(v_i) .
\end{aligned}
\end{equation}
By abuse of notation, $\langle\cdot\rangle$ should be understood as $\textbf{\textit{Q}}(\cdot)$ whenever needed.

For an arbitrary quadrature rule $\textbf{\textit{Q}}$ the $\textbf{\textit{Q}}$-realizable set is defined as \cite{Alldredge2014,Alldredge2015}
$$\textbf{\textit{R}}_N^\textbf{\textit{Q}}=\left\{\gamma~\bigg{|}\gamma=\sum_{i=1}^{n_Q} w_im(v_i)f_i,~f_i>0\right\}.$$
It is worth mentioning that $\textbf{\textit{R}}_N^\textbf{\textit{Q}}$ is a strict (polytopic) subset of $\textbf{\textit{R}}_N$ and like $\textbf{\textit{R}}_N$, it is an open convex cone \cite{Alldredge2014}. $\textbf{\textit{R}}_N^\textbf{\textit{Q}}$ depends on the choice of the quadrature nodes and, in particular, on the number of points $n_Q$.

\subsection{Realizability-preserving property of the schemes}
As we have noted above, the maximum-entropy moment problem \eqref{Momente} has a solution if and only if the moments $\gamma$ are realizable (in $\textbf{\textit{R}}_N$ or $\textbf{\textit{R}}_N^\textbf{\textit{Q}}$, respectively). Since we need to be able to solve the moment problem to evaluate the right hand side of \eqref{ME}, it is crucial to maintain realizability throughout the computation. Similar problems have been treated in the context of radiative transfer, e.g., in \cite{Alldredge2015,Schneider2016a,Schneider2015b,chidyagwai2017comparative,Olbrant2012}.

\begin{theorem}
\label{propm}
Under the non-negativity of the kernels $\Gamma$ and $\omega$ for all $v\in[v_{min},v_{max}]$ and neglected space dependency, the explicit Euler scheme discretization of the moment equations $\eqref{ME}$, combined with $M_N$ or $QMOM_n$, preserves realizability under the CFL-like condition
$$\max\left(\sup_{v\in[v_{min},v_{max}],k}\Big\langle f_{approx}\left(t^k,g_{A^-}(v,v')\right)\omega\left(t^k,v,g_{A^-}(v,v')\right)g_{A^-}'(v,v')\Big\rangle,\sup_{v\in[v_{min},v_{max}],k}\Gamma(t^k,v)\right)\leq\frac{1}{\Delta t},$$
that is, given initial data in $\textbf{\textit{R}}_N^\textbf{\textit{Q}}$ ($\textbf{\textit{R}}_N$, respectively), one step of the forward Euler scheme provides updated values that are also in $\textbf{\textit{R}}_N^\textbf{\textit{Q}}$ ($\textbf{\textit{R}}_N$, respectively).
\end{theorem}
The proof of this theorem can be obtained similarly to the work \cite{Filbet2004, Kumar2006} and \cite{Alldredge2012, Alldredge2015} and is therefore omitted.
\subsection{Implementation details for maximum entropy models}
\label{MaxEn}
For the numerical implementation of the maximum entropy model two big questions arise: how to solve the optimization problem and what to do if the algorithm fails to converge (in a reasonable time). The optimization problem is commonly solved by variations of the Newton algorithm \cite{Alldredge2014, Alldredge2012, Abramov2009,Vie2011}. Problems occur for moment vectors which are near the realizability boundary. In this case, the algorithm may require a large number of Newton iterations to converge or does not converge at all. Also, the solution will be sensitive to small changes in the moments, caused by an ill-conditioned Hessian matrix of the dual objective function. The problem is impaired by the use of a quadrature rule and the finite-precision arithmetic to evaluate the integrals \cite{Alldredge2014}.
There are several variations of Newton's method available that deal with the singularity of the Hessian matrix \cite{Turek1988, Abramowitz1972, Alldredge2012}. We follow the ideas in \cite{Alldredge2014}.  There, the authors use an adaptive change of basis together with a regularization method. 
\subsubsection{Newton's Method with an Adaptive Change of Basis}
\label{NM}
We slightly modify the ideas and algorithms proposed in \cite{Alldredge2014}, where more details and results can be found.
The optimization algorithm we use is based on Newton's method stabilized by an Armijo backtracking line search \cite{Armijo1966} and an adaptive change of basis to improve the condition number of the Hessian.
It computes an approximation $\alpha$ of the true solution $\hat{\alpha}.$ If the algorithm fails to converge we use some regularization technique. 

Let us recall the dual objective function
$$h(\alpha)=\langle G_{\alpha}\rangle-\alpha^T\gamma.$$
Its gradient $g:\mathbb{R}^{N+1}\to\mathbb{R}^{N+1},$ Hessian $H:\mathbb{R}^{N+1}\to\mathbb{R}^{(N+1)\times(N+1)}$ and its Newton direction $d:\mathbb{R}^{N+1}\to\mathbb{R}^{N+1}$ are given by
$$g(\alpha):=\langle mG_{\alpha}\rangle-\gamma,~H(\alpha):=\langle mm^TG_{\alpha}\rangle\text{ and } H(\alpha)d(\alpha)=-g(\alpha).$$
\subsubsection{Initial Guess of the Optimization Algorithm and Regularization} 
First, we define two bases $m_{mono}=(1,v,v^2,...,v^N)^T$ and $m_o$, whereas $m_{mono}$ is used to evaluate the integral terms of \eqref{ME} and to compute the moments. The basis $m_o$ is needed in the optimization algorithm as the  initial orthogonal triangular basis to start the  algorithm. Then the iterative basis remains orthogonal and triangular.
We  choose  as  basis the Legendre basis \eqref{Legendrebasis}. For the Hessian matrix to have full rank it is necessary that $n_Q\geq N+1$ \cite{Alldredge2014}.
In \cite{Alldredge2014} the authors proposed to initially start with the Legendre basis and then save the new adapted basis computed in the Newton algorithm to use it in the next time step as initial basis. Similarly, they handled the approximated Lagrange multipliers $\alpha$.
We proceed  in a slightly different way. We use in every Newton algorithm the Legendre basis $m_o$ as initial basis no matter what time step we have. We then compute an initial guess $\alpha_{M_1}$ for the Lagrange multipliers involving only the first two moments $\gamma_0$, $\gamma_1$ with the help of the Newton algorithm. Here we use as initial guess the Legendre basis and $\alpha_{0,M_1}=\left(\ln(\gamma_0),0\right)^T.$
Our initial guess for the total set of moments is $\alpha_0=\left(\alpha_{M_1},0,...,0\right)^T\in\mathbb{R}^{N+1}.$ Of course for the actual optimization algorithm we need those multipliers with respect to the orthogonal basis $\beta_0=(m_o^{-T}\alpha_0).$ 
If the moments are too close to the realizability boundary, the algorithm could fail, therefore in this case a regularization is used. We define a realizable moment $Q_{mono}$ such that we can derive a new realizable moment which is farther away from the realizability boundary through a convex combination of the original moment $\gamma$ and $Q_{mono}.$ Note that since $\textbf{\textit{R}}_N$, $\textbf{\textit{R}}_N^\textbf{\textit{Q}}$ both are convex cones, convex combination of realizable moments are again realizable.
Therefore, we again take advantage of $\alpha_0$ and compute the corresponding moments  both in the monomial basis and in the Legendre basis
$$Q_{mono}:=\langle m_{mono}\exp\left(\alpha_0^Tm_{mono}\right)\rangle,~~~~~~~Q=\langle m_{o}\exp\left(\beta_0^Tm_{o}\right)\rangle.$$
So starting the Newton algorithm with some moments in the monomial basis $\gamma$ and their corresponding Legendre moments $\tilde{\gamma},$ the regularization technique we use for $r\in[0,1]$ has the following form
$$\tilde{\gamma}_{reg}=(1-r)\tilde{\gamma}+rQ.$$ It is  implemented in this way because for the reduced moment equations \eqref{ME}, mass conservation should be guaranteed, which means $\gamma_1$ remains constant. The Lagrange multiplier $\alpha_0$ exactly reproduces the zeroth, and first moment. So despite the regularization technique they do not change.  
\subsubsection{Stopping Criterion}$\;$ 
For the stopping criterion, we compute the gradient in the $k-$th iteration by
$$g(\alpha_k)=\langle m_{mono}\exp\left(\alpha_k^Tm_{mono}\right)\rangle-\gamma.$$
We stop the algorithm if 
$$\left|\left|\langle m_{mono}\exp\left(\alpha_k^Tm_{mono}\right)\rangle -(1-r)\gamma-rQ_{mono}\right|\right|<\tau.$$
The last two terms arise from the regularization technique. If the algorithm stops and regularization was needed we replace our original moments $\gamma$ by $(1-r)\gamma+rQ_{mono}.$ The rest of the algorithm follows the one shown in \cite{Alldredge2014}.

\subsection{Numerical Treatment of the population balance equation and its moment equations}
We consider the two-dimensional case. Let $\Omega=[x_{min},x_{max}]\times[y_{min},y_{max}]$ be our rectangular spatial domain. We discretize it by a Cartesian equidistant grid with cells $C_{i,j}=[x_{i-\frac{1}{2}},x_{i+\frac{1}{2}}]\times[y_{j-\frac{1}{2}},y_{j+\frac{1}{2}}],$ where $x_{i+1\frac{1}{2}}=x_{min}+i\Delta x,$ $i=0,...,N_x$ and $y_{j+\frac{1}{2}}=y_{min}+j\Delta y,~j=0,...,N_y$ with the step sizes $\Delta x=(x_{max}-x_{min})/N_x$ and $\Delta y$  correspondingly. The cell centers of such a cell $C_{i,j}$ are computed by $(x_i,y_j)=(\frac12(x_{i-\frac{1}{2}}+x_{i+\frac{1}{2}}),\frac12(y_{j+\frac{1}{2}}+y_{j+\frac{1}{2}}))$. The numerical approximation of the average of $f$ over a cell $C_{i,j}$ is defined by
 $$f_{i,j}(t,v)=\frac{1}{\Delta x\Delta y}\int_{C_{i,j}}f(t,x,y,v)dxdy.$$
$T$ is the final time until which the numerical evaluation of the solution is calculated. The time domain $[0,T]$ is discretized by the points of time $t^k=k\Delta t$ with the time step size $\Delta t.$ 
 \subsubsection{Positivity-Preserving Discretization of the Population Balance Equation}
 \label{posPBE}
In this section we derive  a positivity-preserving numerical scheme for the spatially inhomogeneous population balance equations neglecting  here for simplicity the aggregation and breakage terms, which have been treated in a previous section.
 \begin{equation}
 \begin{aligned}
 \label{posFVS}
 \partial_t f+\nabla\cdot(v_df)-\nabla\cdot(D\nabla f)=0.
 \end{aligned}
 \end{equation}
 The droplet velocity is defined by $v_d(t,x)=(u(t,x),z(t,x))^T$ and the diffusion rate by $D=D(t,x)>0.$
 Integrating equation \eqref{posFVS} over the cells $[x_{i-\frac{1}{2}},x_{i+\frac{1}{2}}]\times[y_{j-\frac{1}{2}},y_{j+\frac{1}{2}}]\times[t_k,t_{k+1}]$ and the use of the mid point rule gives the finite volume formulation 
 \begin{equation}
 \begin{aligned}
 \label{disc}
 f_{i,j}^{k+1}(v)=f_{i,j}^k(v)&-\frac{\Delta t}{\Delta x}\bigg[\underbrace{u\left(t^{k},x_{i+\frac{1}{2}},y_j\right)f\left(t^{k},x_{i+\frac{1}{2}},y_j,v\right)}_\text{$=:F_{i+\frac{1}{2},j}^k(v)$}-\underbrace{u\left(t^{k},x_{i-\frac{1}{2}},y_j\right)f\left(t^{k},x_{i-\frac{1}{2}},y_j,v\right)}_\text{$:=F_{i-\frac{1}{2},j}^k(v)$}\bigg]\\&-\frac{\Delta t}{\Delta y}\bigg[\underbrace{z\left(t^{k},x_i,y_{j+\frac{1}{2}}\right)f\left(t^{k},x_i,y_{j+\frac{1}{2}},v\right)}_\text{$:=G_{i,j+\frac{1}{2}}^k(v)$}-\underbrace{z\left(t^{k},x_i,y_{j-\frac{1}{2}}\right)f\left(t^{k},x_i,y_{j-\frac{1}{2}},v\right)}_\text{$:=G_{i,j-\frac{1}{2}}^k(v)$}\bigg]\\
 &+\frac{\Delta t}{\Delta x}\bigg[\underbrace{D\left(t^{k},x_{i+\frac{1}{2}},y_j\right)\partial_xf\left(t^{k},x_{i+\frac{1}{2}},y_j,v\right)}_\text{$:=H_{i+\frac{1}{2},j}^k(v)$}-\underbrace{D\left(t^{k},x_{i-\frac{1}{2}},y_j\right)\partial_xf\left(t^{k},x_{i-\frac{1}{2}},y_j,v\right)}_\text{$:=H_{i-\frac{1}{2},j}^k(v)$}\bigg]
 \\&+\frac{\Delta t}{\Delta y}\bigg[\underbrace{D\left(t^{k},x_i,y_{j+\frac{1}{2}}\right)\partial_yf\left(t^{k},x_i,y_{j+\frac{1}{2}},v\right)}_\text{$:=M_{i,j+\frac{1}{2}}^k(v)$}-\underbrace{D\left(t^{k},x_i,y_{j-\frac{1}{2}}\right)\partial_yf\left(t^{k},x_i,y_{j-\frac{1}{2}},v\right)}_\text{$:=M_{i,j-\frac{1}{2}}^k(v)$}\bigg].
  \end{aligned}
  \end{equation}
We exemplary show for $F_{i+\frac{1}{2},j}^k$ and $H_{i+\frac{1}{2},j}^k$  how to handle these terms. $G$ and $M$ work similarly.
Let us first address the advection flux in $x$-direction $F_{i+\frac{1}{2},j}^k$. We split droplets into left- and right-moving particles and get for the numerical flux
$$F_{i+\frac{1}{2},j}^k(v)=\max\left(0,u_{i+\frac{1}{2},j}^k\right)f^{+,k}_{i+\frac{1}{2},j}(v)+\min\left(0,u_{i+\frac{1}{2},j}^{k}\right)f^{-,k}_{i+\frac{1}{2},j}(v),$$
where $u_{i+\frac{1}{2},j}^k=1/2\left(u_{i,j}^k+u_{i+1,j}^k\right).$ The values $f^{+,k}_{i+\frac{1}{2},j}$, $f^{-,k}_{i+\frac{1}{2},j}$ on the right and left sides of the cell edge at $(x_{i+\frac{1}{2}},y_j)$ \cite{Hauck2011, Garrett2013}, are approximated by
$$f_{i+\frac{1}{2},j}^{+,k}(v)=f_{i,j}^k(v)+\frac{1}{2}\sigma^{x,k}_{i,j}(v),~~~f_{i+\frac{1}{2},j}^{-,k}(v)=f_{i+1,j}^k(v)-\frac{1}{2}\sigma_{i+1,j}^{x,k}(v),$$
with $\sigma_{i,j}^{x,k}$ the approximation of the slope in the $x-$direction in cell $C_{i,j}$ at time $t_k,$
$$\sigma_{i,j}^{x,k}(v)=\text{minmod}\left\{2\left(f_{i+1,j}^k(v)-f_{i,j}^k(v)\right),\frac{1}{2}\left(f_{i+1,j}^k(v)-f_{i-1,j}^k(v)\right),2\left(f_{i,j}^k(v)-f_{i-1,j}^k(v)\right)\right\}.$$
The minmod limiter does not only guarantee non-negativity of $f^{+,k}_{i+\frac{1}{2},j}$, $f^{-,k}_{i+\frac{1}{2},j}$ but also suppresses spurious oscillations \cite{Toro}.
The diffusion term is treated in a similar way as in \cite{Chertock2008},
$$H_{i+\frac{1}{2},j}^k(v)=\frac{1}{2\Delta x}\left(D_{i+1,j}^{k+1}+D_{i,j}^k\right)\left(f_{i+1,j}^k(v)-f_{i,j}^k(v)\right).$$
\subsubsection{Realizability Preserving Discretization of the Moment Equations}
For the maximum entropy closure or the quadrature method of moments we did discuss the realizability preservation of the spatially homogeneous case in subsection \ref{Mereal}. Here, a  realizability preserving scheme for the advection-diffusion equation neglecting  the breakage and aggregation terms
is given. Consider
\begin{equation}
\begin{aligned}
\label{posME}
\partial_t \gamma+\nabla(v_d\gamma)-\nabla(D\nabla\gamma)=0.
\end{aligned}
\end{equation}
 Assuming that $\gamma_{i,j}^k$ is realizable and that $\tilde{f}$ is the corresponding distribution function approximated by the maximum entropy model. Then to derive a realizable scheme we multiply equation \eqref{disc} with the basis $m$ and integrate over the whole volume domain \cite{Hauck2011,Schneider2015b}, which gives
  \begin{equation}
 \begin{aligned}
 \label{schmememax}
 \gamma_{i,j}^{k+1}= \gamma_{i,j}^k&-\frac{\Delta t}{\Delta x}\left\langle\tilde{F}^{k}_{i+\frac{1}{2},j}(v)-\tilde{F}^{k}_{i-\frac{1}{2},j}(v)\right\rangle-\frac{\Delta t}{\Delta y}\left\langle\tilde{G}^{k}_{i,j+\frac{1}{2}}(v)-\tilde{G}^{k}_{i,j-\frac{1}{2}(v)}\right\rangle\\
&+\frac{\Delta t}{\Delta x}\left\langle\tilde{H}_{i+\frac{1}{2},j}^k(v)-\tilde{H}_{i-\frac{1}{2},j}^k(v)\right\rangle+\frac{\Delta t}{\Delta y}\left\langle\tilde{M}_{i,j+\frac{1}{2}}^k(v)-\tilde{M}_{i,j-\frac{1}{2}}^k(v)\right\rangle,
 \end{aligned}
 \end{equation}
 the terms $\tilde{F}_{i\pm\frac{1}{2},j}^{k},~\tilde{G}_{i,j\pm\frac{1}{2}}^{k},~\tilde{H}_{i\pm\frac{1}{2},j}^{k},~\tilde{M}_{i,j\pm\frac{1}{2}}^{k}$ are in the same way defined as $F_{i\pm\frac{1}{2},j}^{k},~G_{i,j\pm\frac{1}{2}}^{k},~H_{i\pm\frac{1}{2},j}^{k},~M_{i,j\pm\frac{1}{2}}^{k}$ but with respect to the approximate distribution function $\tilde{f}$.  
For the polynomial closure $P_N$ we can also use this scheme, but there no realizability can be guaranteed, since the distribution function can get negative. \\
A  similar scheme as in  \eqref{schmememax}
applied to  the  moment equations \eqref{posME} with the $QMOM$ closure without  aggregation and breakage terms can be shown to be  realizability preserving in every time step  \cite{Vikas2010, Vikas2013}. Recall that the $QMOM-$distribution function is of the form $\tilde{f}(t,x,v)=\sum_{\alpha}w_\alpha(t,x)\delta\left(v-v_\alpha(t,x)\right).$
The discretization is now similar to the discretization method \eqref{schmememax}. One only has to change the definition of the terms $\tilde{f}_{i\pm\frac{1}{2},j}^{+/-,k},~\tilde{f}_{i,j\pm\frac{1}{2}}^{+/-,k}$.
They are determined by
$$\tilde{f}_{i+\frac{1}{2},j}^{+,k}(v)=\sum_{\alpha}\tilde{w}_{\alpha,i+\frac{1}{2},j}^{+,k}\delta\left(v-v_{\alpha,i,j}^k\right),~\tilde{f}_{i+\frac{1}{2},j}^{-,k}(v)=\sum_{\alpha}\tilde{w}_{\alpha,i+\frac{1}{2},j}^{-,k}\delta\left(v-v_{\alpha,i+1,j}^k\right).$$
The interface weights are defined as
$$\tilde{w}_{\alpha,i+\frac{1}{2},j}^{+,k}=w_{\alpha,i,j}^{k}+\frac{1}{2}\partial w_{\alpha,i,j}^{x,k},~\tilde{w}_{\alpha,i+\frac{1}{2},j}^{-,k}=w_{\alpha,i+1,j}^{k}-\frac{1}{2}\partial w_{\alpha,i+1,j}^{x,k}.$$
Again we use the minmod limiter to reconstruct the slopes as
$$\partial w_{\alpha,i,j}^{x,k}=\text{minmod}\left(2\left(w_{\alpha,i+1,j}^{k}-w_{\alpha,i,j}^{k}\right), \frac{1}{2}\left(w_{\alpha,i+1,j}^{k}-w_{\alpha,i-1,j}^{k}\right),2\left(w_{\alpha,i,j}^{k}-w_{\alpha,i-1,j}^{k}\right)\right).$$ 
The discretization in the $y-$direction works analogously.
  \begin{theorem}
To ensure either a non-negative distribution function for the advection-diffusion equation \eqref{posFVS} discretized by the finite volume scheme \eqref{disc} or a realizability preservation of the schemes for the moment equations with $QMOM$, or the maximum entropy closures of the advection diffusion problem \eqref{posME}, the following time step restriction has to be fulfilled
$$\sup_{\theta\in(0,1)}\left(\frac{2a}{\theta\Delta x},\frac{2b}{\theta\Delta y},\frac{1}{2\left(1-\theta\right)}\left(\frac{c}{\Delta x^2}+\frac{d}{\Delta y^2}\right)\right)\leq\frac{1}{\Delta t},$$
with 
$$a=\sup_{i,j,k}\left(\left|u_{i+1,j}^k+u_{i,j}^k\right|\right),~b=\sup_{i,j,k}\left(\left|z_{i,j+1}^k+z_{i,j}^k\right|\right),$$
and $$c=\sup_{i,j,k}\left(D_{i+1,j}^k+2D_{i,j}^k+D_{i-1,j}^k\right),~d=\sup_{i,j,k}\left(D_{i,j+1}^k+2D_{i,j}^k+D_{i,j-1}^k\right).$$
\end{theorem}

%% file: Sections/results.tex
\section{Numerical Results}
\label{results}
We compare the zeroth moment of the different closures $QMOM$ \eqref{QMOMf}, maximum entropy \eqref{max2} and $P_N$ \eqref{PN} with each other. As reference solution we use the finite volume scheme explained in section \ref{FVSchap}. The zeroth moment for the finite volume scheme is computed by
$$\gamma_0=\Delta v\sum_{i=1}^{n_v}f_i.$$
We will compare them regarding the $L_2$-error as well as computation time.
The relative error in the $L^2$ sense for two number densities $\gamma_{0,\text{ref}}$, $\gamma_{0,\text{clo}}$ is evaluated by the formula
\begin{align}
\label{eq:relL2}
E_2(\gamma_{0,\text{ref}}, \gamma_{0,\text{clo}})^2=\frac{\int_{0}^T\int_{\Omega}\left|\gamma_{0,\text{ref}}(t,x)-\gamma_{0,\text{clo}}(t,x)\right|^2dxdt}{\int_{0}^T\int_{\Omega}\left|\gamma_{0,reference}(t,x)\right|^2dxdt},
\end{align}
where $\Omega$ is the spatial domain and $T$ the final time until which we compute the numerical solution. 
For the Newton-algorithm in the maximum entropy case we need to choose several parameters which can be found in table \ref{tablemax}. For the explanation of the parameters see \cite{Alldredge2014}.
\begin{table}[!htb]
\begin{center}
\begin{tabular}{|c|ccccc|}
\hline
Parameter & $k_{max}$ &$\epsilon$ &$\chi$ &$\tau$ &$\{r_l\}$ \\
\hline
Value &$400$& $2^{-52}$ &$3/5$ & $10^{-9}$& $\{0,10^{-10},10^{-8},10^{-6},10^{-4},10^{-2},0.1,0.5,1\}$  \\
\hline
\end{tabular}
\caption{Parameter choice for the Newton algorithm in the case of the maximum entropy closure.}
\label{tablemax}
\end{center}
\end{table}

For all examples we use the a daughter droplet distribution which is based on the purely statistical daughter droplet distribution function of Hill and Ng \cite{Hill2004}. Recall that $N(v')$ is the average number a particle of size $v'$ splits. Our daughter droplet distribution is a weighted sum of daughter droplet functions
$$\beta(v,v')=\left\{\begin{array}{cl} 0 & \mbox{for $v<v_{min},v'< v_{min}$,}\\ \sum_{i=1}^{c(v')}g_i(v')\bar{\beta}_i(v,v') &\mbox{for $v\in[v_{min},v']$ and $v'\in[v_{min},v_{max}]$} \\ 0& \mbox{for $v'>v_{max}$}\end{array}\right.,$$
where the $i-$breakup process (a mother droplet splits in average into $i\in\mathbb{N}$ smaller droplets) daughter distribution functions for $v_{max}\geq v'> pv_{min}$ are defined as 
\begin{align*}
\bar{\beta}_i(v,v')=\left\{\begin{array}{cl} 
0& \mbox{for $v>v'-(i-1)v_{min}$}\\
\tilde{ \beta}_i(v,v')  & \mbox{for }v_{min}\leq v\leq v'-(i-1)v_{min}
 \end{array}\right.,
\end{align*}
\begin{align*}
\label{ddd}
\text{with }\hspace{0.5cm}\tilde{\beta}_i(v,v')=\frac{i(mi+i-1)!(v-v_{min})^m(v'-v-(i-1)v_{min})^{mi+i-m-2}}{m!(mi+i-m-2)!(v'-iv_{min})^{im+i-1}},
\end{align*}
where $m\in\mathbb{N}$ determines the shape of $\tilde{\beta}_i.$ The no-breakup daughter distribution function is defined for $v_{max}\geq v'\geq v_{min}$ as
\begin{align*}
\bar{\beta}_1(v,v')=\left\{\begin{array}{cl} 
0& \mbox{for $v>v'$}\\
 \delta(v-v')  & \mbox{for }v_{min}\leq v\leq v'
 \end{array}\right..
\end{align*}
For all $v'\in[v_{min},v_{max}]$ the weight functions $g_i(v')\in[0,1]$ have to fulfill 
\begin{equation}
\begin{aligned}
\label{dropletbreak}
\sum_{i=1}^{c(v')}g_i(v')=1\text{, }\sum_{i=1}^{c(v')}g_i(v')i=N(v')\text{ and } g_i(v')=0\text{ for }v'<iv_{min},
\end{aligned}
\end{equation}
 and $c(v')\in[2,\infty)$. 
We choose $c(v')$ as $c(v')=2N(v')-1.$ For the computation of the weight functions let $\textbf{1}_U(x)$ be the indicator function on $U$. We choose $g_{N(v')-j}=g_{N(v')+j}$ for all $j=1,...,N(v')-1$ and for all $v'\in[v_{min},v_{max}]$ and assume that the weight functions are of the form
$$g_{i}(v')=\textbf{1}_{]0,\infty[}(i+1-k(v',N(v')))\frac{1}{j2^{N(v')-i}},\text{ for }i=1,...,N(v'),$$
with $k(v',N(v'))\in\mathbb{N}$ determining that the first $k(v',N(v'))-1$ and last $k(v',N(v'))-1$ weight functions of the series $g_1(v'),..,g_{2N(v')-1}(v')$ are zero for $v'$.
With the first property \eqref{dropletbreak} of the weight functions, the variable $j$ is determined by
$$j=3-\frac{1}{2^{N(v')-1-k(v',N(v'))}}.$$
Let us define the intervals 
$$I_{1}=[v_{min},2v_{min}]\text{, } I_{l}=]lv_{min},(l+1)v_{min}],~ l=2,...,2p-2, \text{ and }I_{2p-1}=](2p-1)v_{min},v_{max}].$$
Then, if $v'\in I_l$ for some $l\in\mathbb{N}$ we define
$$k(v',N(v')):=\left\{\begin{array}{cl} l & \mbox{for $l\in[1,...,p]$,}\\ 2p-l& \mbox{for $l\in[p+1,...,2p-2]$}\\ 1& \mbox{for $l=2p-1$}\end{array}\right.,N(v'):=\left\{\begin{array}{cl} l & \mbox{for $l\in[1,...,p]$,}\\ p& \mbox{for $l\in[p+1,...,2p-2]$}\\ p& \mbox{for $l=2p-1$}\end{array}\right..$$
This definition  fulfills all the properties mentioned in section \ref{PBEsec}~.

\begin{figure}
\begin{center}
\subfloat[][$v'\in I_1$]{
	\includegraphics[width=0.2\textwidth]{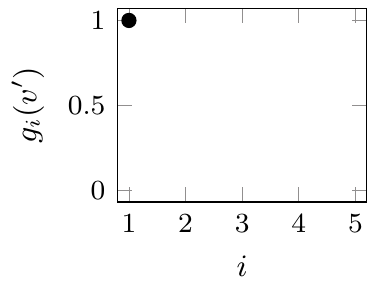}
}	
\subfloat[][$v'\in I_2$]{
\includegraphics[width=0.15\textwidth]{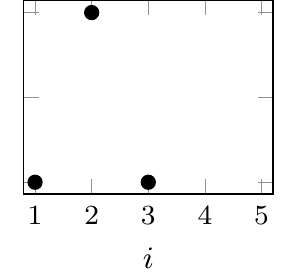}
}
\subfloat[][$v'\in I_3$]{
\includegraphics[width=0.15\textwidth]{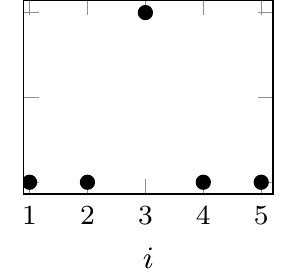}
}
\subfloat[][$v'\in I_4$]{
\includegraphics[width=0.15\textwidth]{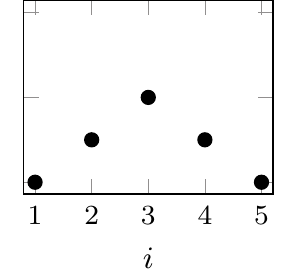}
}
\subfloat[][$v'\in I_5$]{
\includegraphics[width=0.15\textwidth]{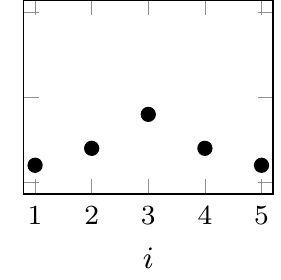}}
\caption{
The subfigures (a)-(e) show the distribution of the weights $g_i(v')$ for quasi-ternary breakage (p=3) and for $v'$ contained in the different intervals $I_1$-$I_5$. It can be seen that only for $v'\in I_j$, $j>3$ different kind of breakages occur at the same time. }
\label{figure dropletbreak}
\end{center}
\end{figure}

Although the choice of the weight functions seems complicated there are just chosen such that they are symmetric about $N(v')$ and the largest weight corresponds to $N(v')$ for all $v'$ and such that all properties mentioned in section \ref{PBEsec} are fulfilled. The choice also guarantees that droplets cannot break into droplets whose volume is not contained in $\left[v_{min},v_{max}\right]$. Figure \ref{figure dropletbreak} demonstrates the behavior of the weight functions $g_i(v')$ on the example of quasi-ternary breakage.
\noindent

\subsection{Pure Breakage}
\begin{figure}

		\begin{center}
				\subfloat[][]{\includegraphics[width=0.49\textwidth]{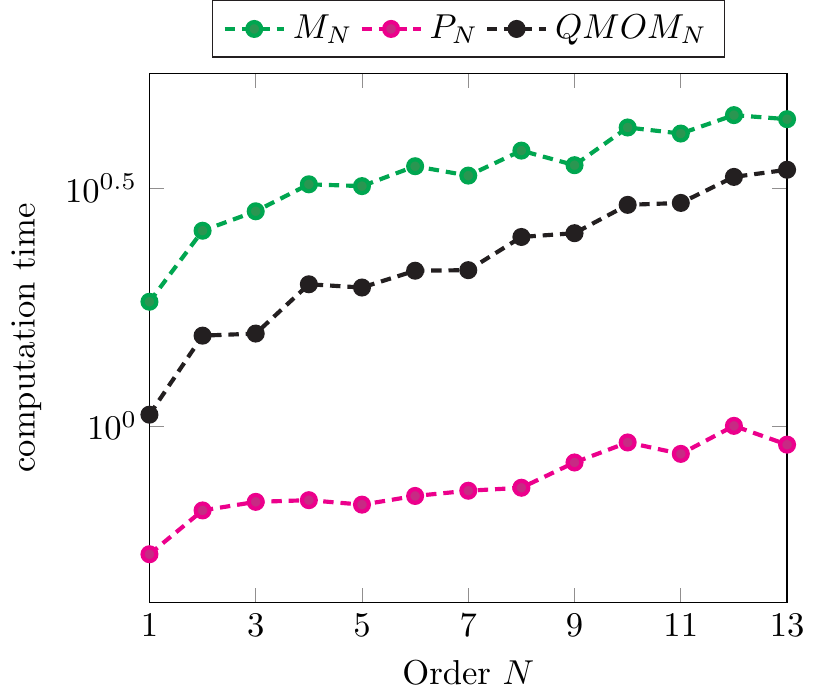}}
		\subfloat[][]{\includegraphics[width=0.49\textwidth]{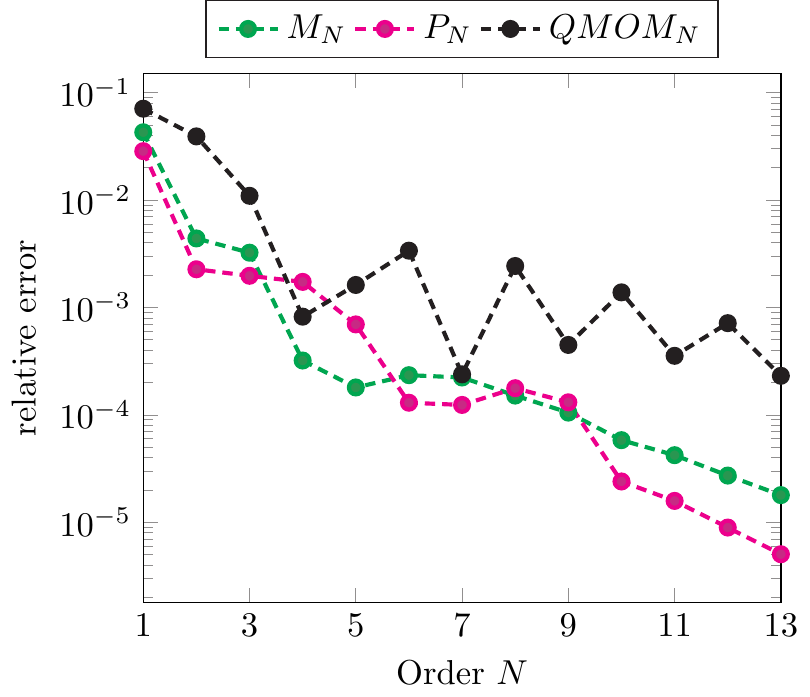}}
	
	\subfloat[][]{\includegraphics[width=0.49\textwidth]{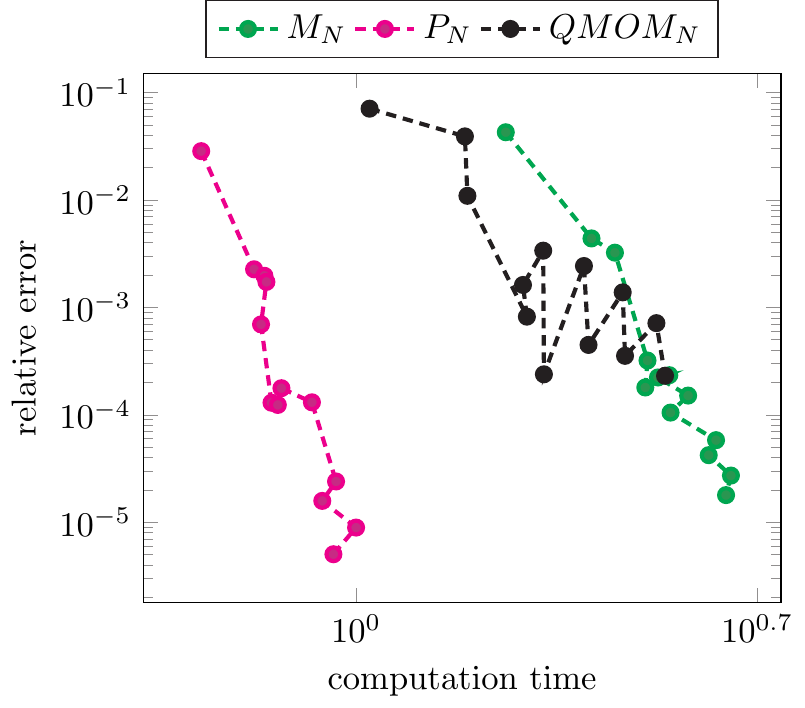}}	
	\caption{Pure breakage: In subfigure (a) the order of moments $N$ is plotted against the computation time, for the maximum entropy model, the polynomial closure and the $QMOM$ and in subfigure (b) the order $N$ is plotted against the relative $E_2$ error for all three methods, in subfigure (c) we can see the computation time plotted against the relative $E_2$ error, again for all three models. In all three subfigures the computations are only done up to the thirteenth order. The reference solution for the error evaluation is computed by the finite volume scheme (FVS) with $n_v=5000.$}
			\label{figure break}
					\end{center}
		\end{figure}

First we only want to consider a pure breakage problem without spatial dependency. This means that we set $\omega\equiv0$ in \eqref{ME} and \eqref{PBE}.
Our initial condition is a normal distribution
\begin{equation}
\begin{aligned}
\label{init}
f(0,v)=N_0/\left(\sqrt{2\pi}\sigma_{init}\right)\exp\left(-1/2\left(v-v_0\right)^2/\sigma_{init}^2\right).
\end{aligned}
\end{equation}
For the daughter droplet distribution function we choose binary breakage $p=2$ and $m=2$.
Table \ref{tablebreak} contains some of the parameters chosen for this example.
\begin{table}[!htb]
\begin{center}
\begin{tabular}{|c|ccccccc|}
\hline
Parameter & $v_{min}$ &$v_{max}$ & $\sigma_{init}$ &$\alpha_0$ &$d_0$ &$N_0$ &$v_0$\\
\hline
Value &$\frac{1}{6}\pi0.001^3$& $\frac{1}{6}\pi0.009^3$ &$0.1(v_{max}-v_{min})$ & $0.01$& $\left(\frac{3}{\pi}(v_{min}+v_{max})\right)^{(1/3)}$ &$\frac{6\alpha_0}{\pi d_0^3}$& $\frac{\pi}{6}d_0^3$ \\
\hline
\end{tabular}
\caption{Parameter choice for the volume domain and the initial condition.}
\label{tablebreak}
\end{center}
\end{table}\\
The time step size is for all methods $\Delta t=0.01$. We compute the solutions until $T=4$ is reached. The number of points of the quadrature rule is $n_Q=100$ and for the reference solution computed by the finite volume scheme we use $n_v=5000$ grid points.
As breakage frequency we choose the one from Coulaloglou and Tavlarides \cite{Coulaloglou1977}
\begin{equation}
\begin{aligned}
\label{freq}
\Gamma(v)=\frac{C_1\epsilon^{\frac{1}{3}}}{\left(1+\alpha_d\right)v^{\frac{2}{9}}}\exp\left(-\frac{C_2\sigma\left(1+\alpha_d\right)^2}{\rho_d\epsilon^{\frac{2}{3}}v^{\frac{5}{9}}}\right).
\end{aligned}
\end{equation}
Table \ref{tablegamma} shows the values corresponding to the parameters for the breakage frequency.
 \begin{table}[!htb]
\begin{center}
\begin{tabular}{|c|cccccc|}
\hline
Parameter &$C_1$ & $C_2$ & $\alpha_d$ &$\epsilon$ &$\rho_d$ & $\sigma$\\
\hline
Value & $0.12$ &$0.078$&$\alpha_0$ & $0.004$&  $865.6$ & $0.0361$ \\
\hline
\end{tabular}
\caption{Parameter choice for the breakage frequency.}
\label{tablegamma}
\end{center}
\end{table}
In Figure \ref{figure break} the results are shown for the comparison of $P_N$ (magenta line), $M_N$ (green line) and $QMOM_N$ (black line). We computed the density for all three models with different values for the order, namely up to order thirteenth. The figure is divided into three subfigures: (a) showing the order $N$ plotted against the computation time, (b) illustrating the relation between the order and the relative $L_2$ error \eqref{eq:relL2} and (c) showing the computation time plotted against the relative $L_2$ error \eqref{eq:relL2}. It is demonstrated that with increasing order the computation time also increases. For this situation $M_N$ and  $QMOM_N$
have much higher computation times than  $P_N$. 
We note that the  optimization algorithm for  $M_N$ is the dominating part of 
the computation. In case of $QMOM_N$ the solution of the nonlinear system is the costly part. This results  in the higher computation times for $M_N$ and $QMOM_N$. On the other hand, as mentioned above, the 
$P_N$ approach  might lead to unphysical moments. 

		\noindent
 \subsection{Pure Aggregation}
\begin{figure}

		\begin{center}
				\subfloat[][]{\includegraphics[width=0.49\textwidth]{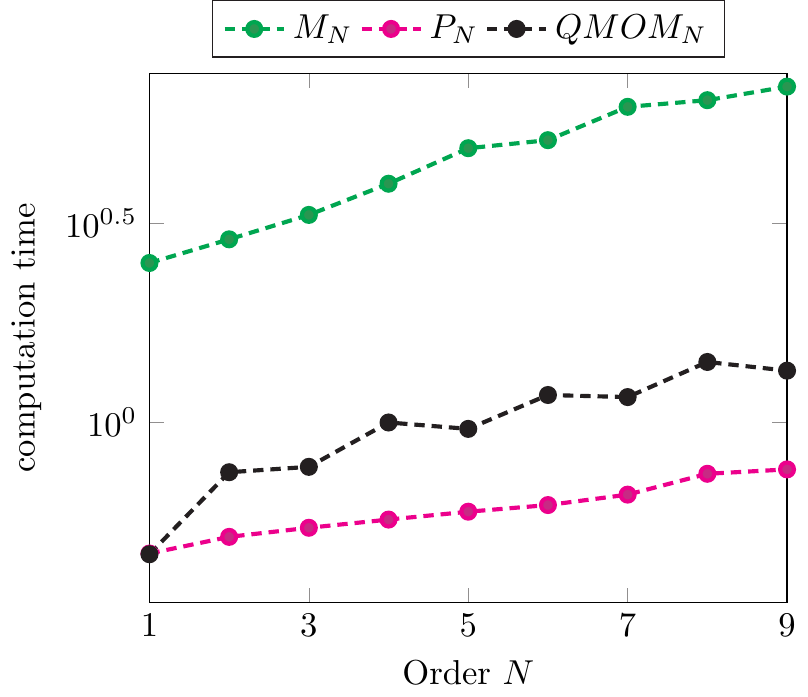}}
		\subfloat[][]{\includegraphics[width=0.49\textwidth]{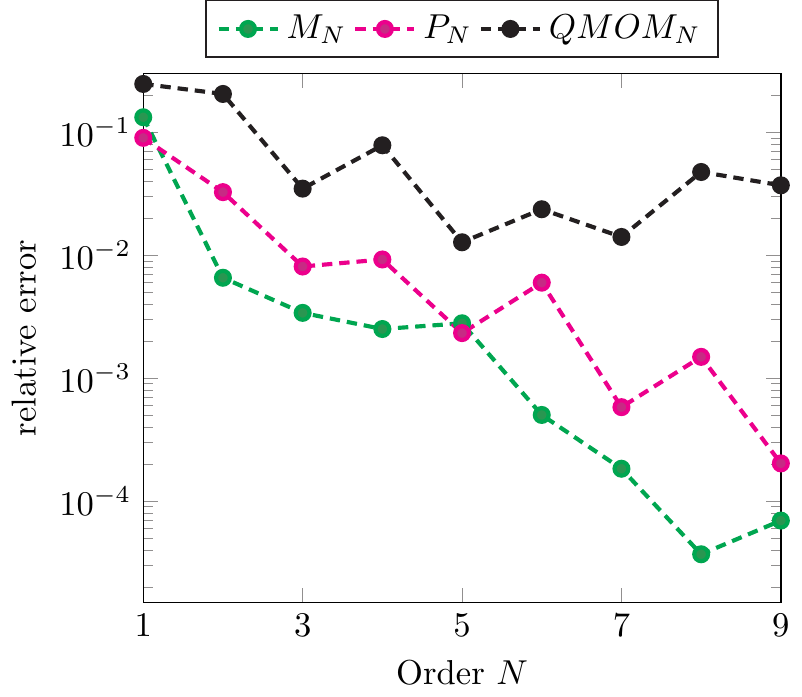}}
	
	\subfloat[][]{\includegraphics[width=0.49\textwidth]{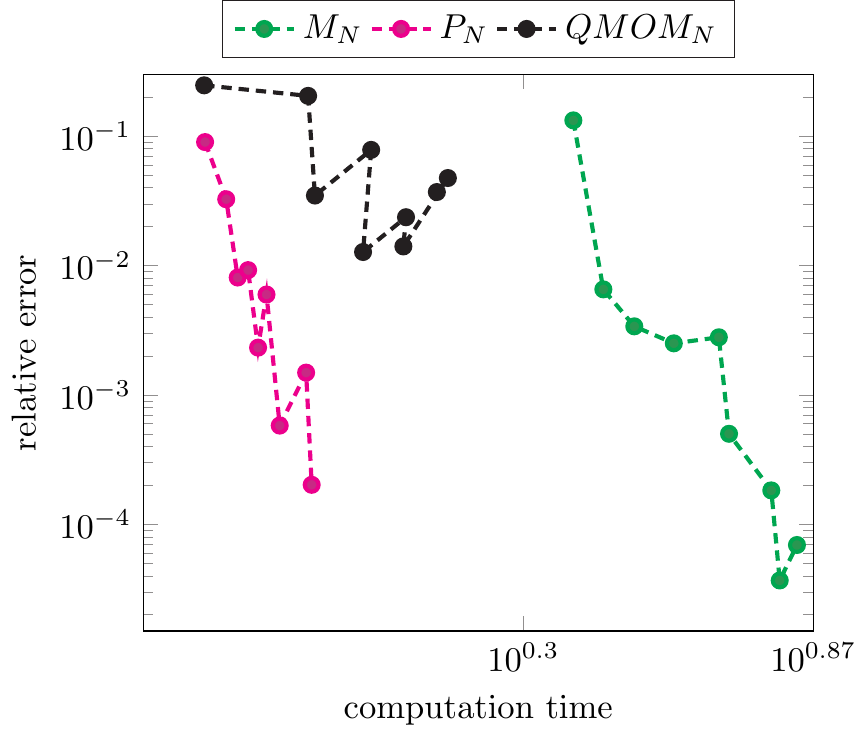}}	
	\caption{Pure aggregation: In subfigure (a) the order of moments $N$ is plotted against the computation time, for the maximum entropy model, the polynomial closure and the $QMOM$ and in subfigure (b) the order $N$ is plotted against the relative $E_2$ error for all three methods, in subfigure (c) we can see the computation time plotted against the relative $E_2$ error, again for all three models. In all three subfigures the computations are  only done up to the ninth order. The reference solution for the error evaluation is computed by the finite volume scheme (FVS) with $n_v=5000.$}
			\label{figure agg}
					\end{center}
		\end{figure}

For the pure aggregation problem we choose $\Gamma\equiv0,~\beta\equiv0$ in \eqref{ME} and \eqref{PBE} and neglect spatial dependency. We use an aggregation kernel proposed by Coulaloglou and Tavlarides \cite{Coulaloglou1977}
\begin{equation}
\begin{aligned}
\label{aggkernel}
 \omega(v,w)=\frac{C_{\Omega}}{1+\alpha_d}\left(v^{\frac{1}{3}}+w^{\frac{1}{3}}\right)^2\epsilon^{\frac{1}{3}}\left(v^{\frac{2}{9}}+w^{\frac{2}{9}}\right)^{\frac{1}{2}}\exp\left(-\frac{k_{\omega'}\eta_c\rho_c\epsilon}{\sigma^2\left(1+\alpha_d\right)^3}\left(\frac{v^{\frac{1}{3}}w^{\frac{1}{3}}}{v^{\frac{1}{3}}+w^{\frac{1}{3}}}\right)^4\right).
 \end{aligned}
 \end{equation}
The parameters remain the same as in the pure breakage case see table \ref{tablebreak}. Table \ref{tableomega} shows the parameters for the aggregation kernel.
 \begin{table}[!htb]
\begin{center}
\begin{tabular}{|c|ccccccc|}
\hline
Parameter &$C_{\Omega}$ &$k_{\omega'}$ &$\alpha_d$  &$\epsilon$ &$\rho_c$ & $\sigma$ &$\eta_c$\\
\hline
Value &$41.2$& $1.33\cdot10^{10}$ &$\alpha_0$ & $0.004$&  $1000$ & $0.0361$ &$0.001$ \\
\hline
\end{tabular}
\caption{Parameter choice for the aggregation kernel.}
\label{tableomega}
\end{center}
\end{table}\\
As before, we choose $n_Q=100$, $\Delta t=0.01$, $n_v=5000$.
The comparison of $P_N$ (magenta line), $M_N$ (green line) and $QMOM_N$ (black line) is illustrated in figure \ref{figure agg}. Again the densities for all three models are computed with changing values for the order, namely up to the ninth order. The figure is subdivided in the same way as figure \ref{figure break}. It is demonstrated that with increasing order the computation time also increases. 
In this case again the computation times of $M_N$ and $QMOM_N$ are much larger than those of $P_N$.
As before the optimization algorithm is the dominating part of $M_N$, resulting in the higher computation time. 
 If we compare  figures \ref{figure break} and \ref{figure agg}, one observes that the relative error in the $QMOM_N$ case is much higher. This is due to the fact, that the breakage integrals are less complex then the aggregation integrals, where the distribution function has to be evaluated twice.
		\noindent
		 \subsection{Coupled Diffusion-Advection-Breakage-Aggregation} 
\begin{figure}	
	\begin{center}
		\captionsetup[subfloat]{farskip=0.01pt,captionskip=0.10pt}	
\includegraphics[width=0.9\textwidth]{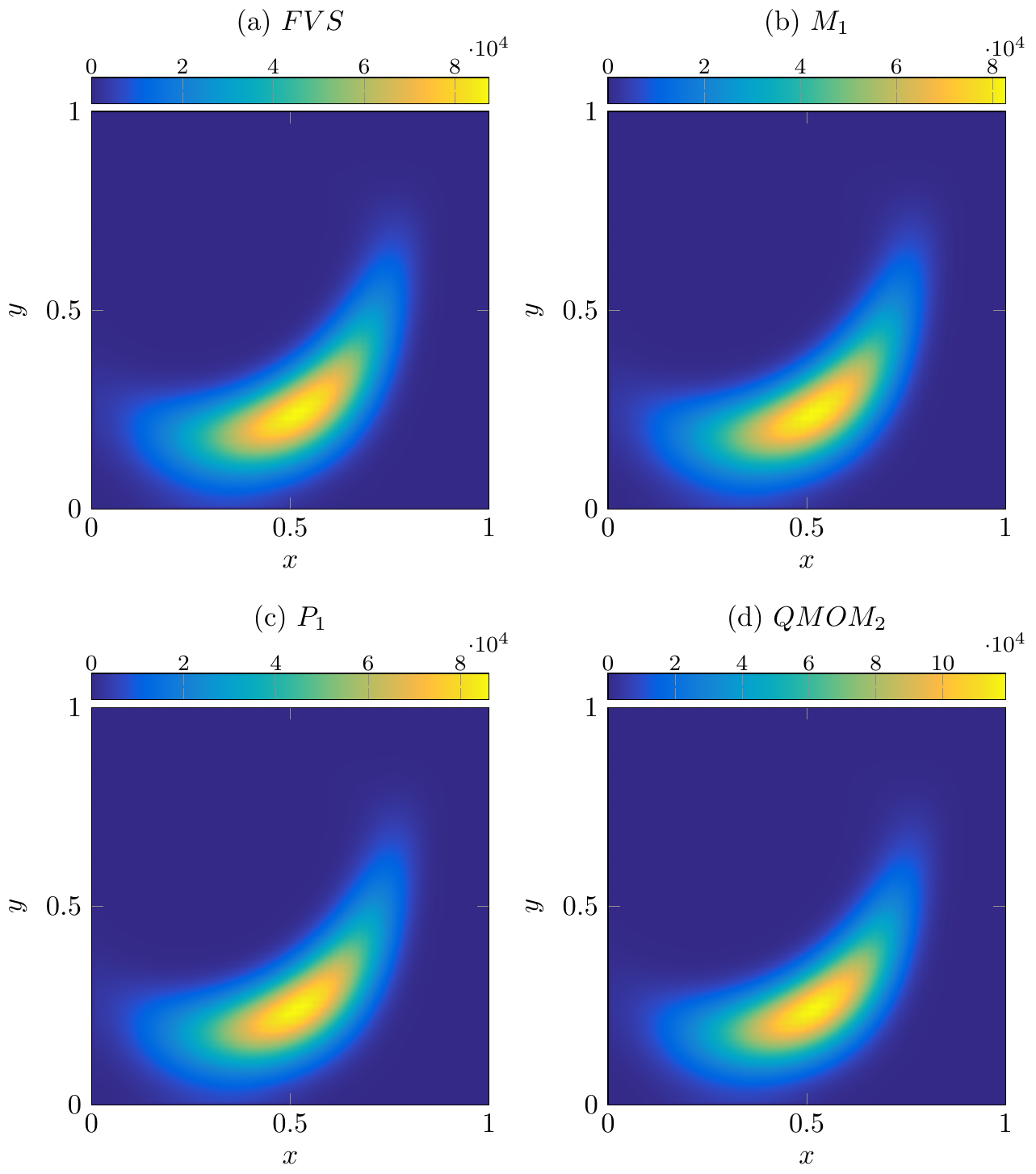}
\end{center}
			\caption{ Densities at the point of time $t=2.1$.}	
				\label{figure1}
		\end{figure}

\begin{figure}	
	\begin{center}
		\captionsetup[subfloat]{farskip=0.01pt,captionskip=0.10pt}	
\includegraphics[width=0.9\textwidth]{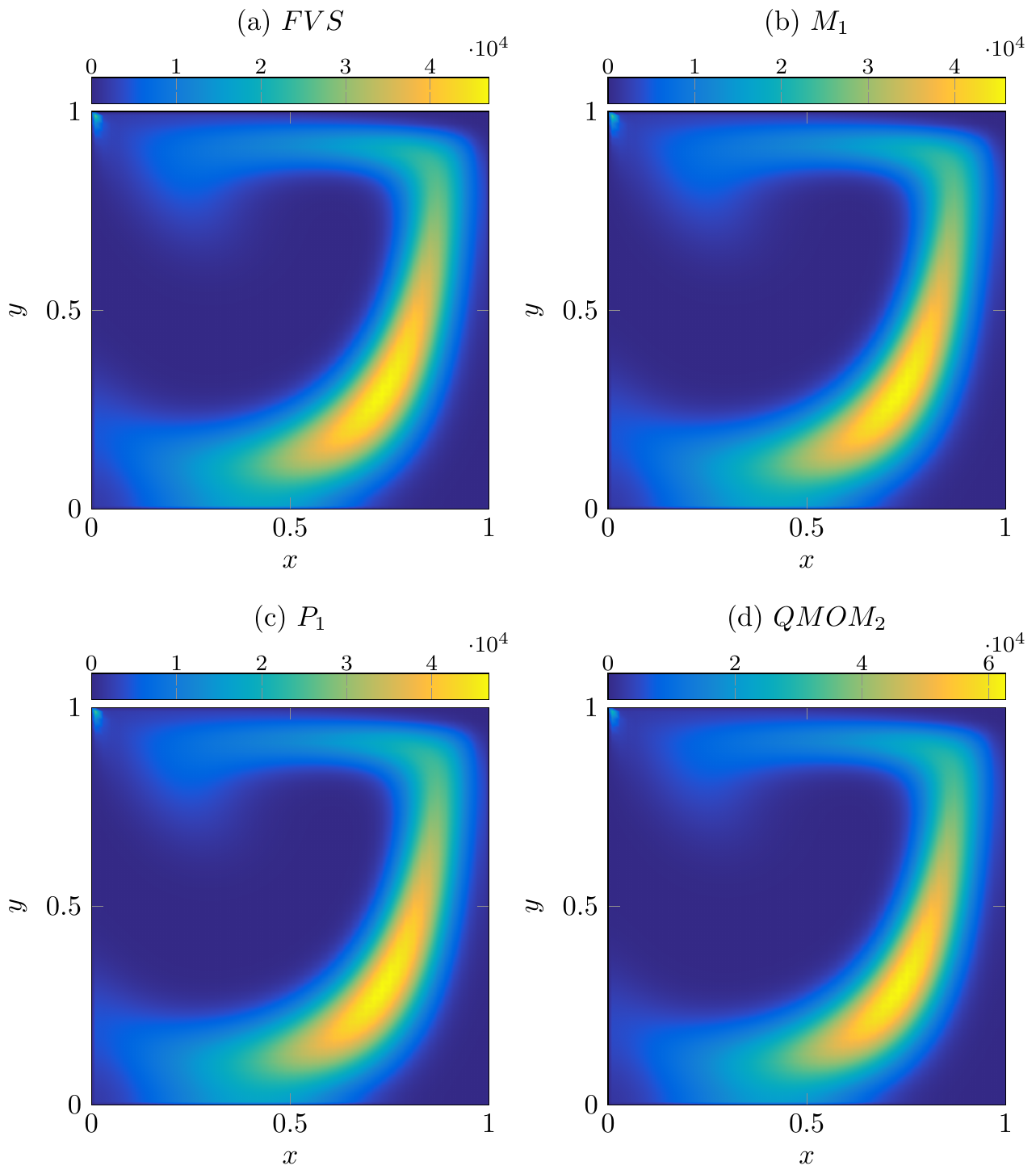}

		\end{center}
			\caption{ Densities at the point of time $t=4.6$.}	
				\label{figure2}
		\end{figure}

For this example we consider equations \eqref{PBE} and \eqref{ME}, we choose the diffusion rate to be constant $D(t,x)=D=0.001.$ We assume that the droplet velocity $v_d(t,x)$ is the same as the velocity of the fluid computed by the lid driven cavity problem \cite{Bruneau2006a} with a Reynold's number of $Re=5.$
 The spatial domain is a square $\Omega=[0,1]\times[0,1]$. We discretize it with $50$ gridpoints in both directions. We choose as breakage frequency \eqref{freq} with the same parameters as in table \ref{tablegamma}. For the daughter droplet distribution function we choose again $p=2$ and $m=2$. The aggregation kernel is the one from \eqref{aggkernel} of the pure aggregation problem with the parameters from table \ref{tableomega}. Again, $n_Q=100.$ The reference solution is computed via the finite volume scheme from section \ref{FVSchap} with $n_v=2000$.\\
 We compare the densities of the finite volume schemes with those of the moment methods, namely $QMOM$ \eqref{QMOMf}, maximum entropy \eqref{max2} and $P_N$ \eqref{PN}. The numerical simulation for all methods runs until the final time $T=5$ is reached.
As initial distribution function we choose a Gaussian distribution, 
$$f(0,x,y,v)=\frac{1}{4\pi^20.08^2}\exp\left(-\frac{1}{2}\left(\frac{(x-0.3)^2}{0.08^2}+\frac{(y-0.3)^2}{0.08^2}\right)\right)\frac{N_0}{\sqrt{2\pi}\sigma_{init}}\exp\left(-1/2\frac{\left(v-v_0\right)^2}{\sigma_{init}^2}\right),$$
the parameters for $f$ are chosen as in table \ref{tablebreak}.\\
In figure \ref{figure1} and \ref{figure2} the time evolution of the densities for the different methods, namely FVS (a), $M_1$ (b), $P_1$ (c), $QMOM_2$ (d) are exemplary shown for two different points of time.\\
With increasing time the densities move from the lower left corner to the upper left corner, following the surrounding fluid. We can see that $P_1$ and $M_1$ look similar to the reference solution computed by the FVS. For $QMOM_2$ the solution shows bigger differences compared to  the reference solution. 

	\setcounter{subfigure}{0} 

		\noindent

\begin{figure}

		\begin{center}
				\subfloat[][]{\includegraphics[width=0.45\textwidth]{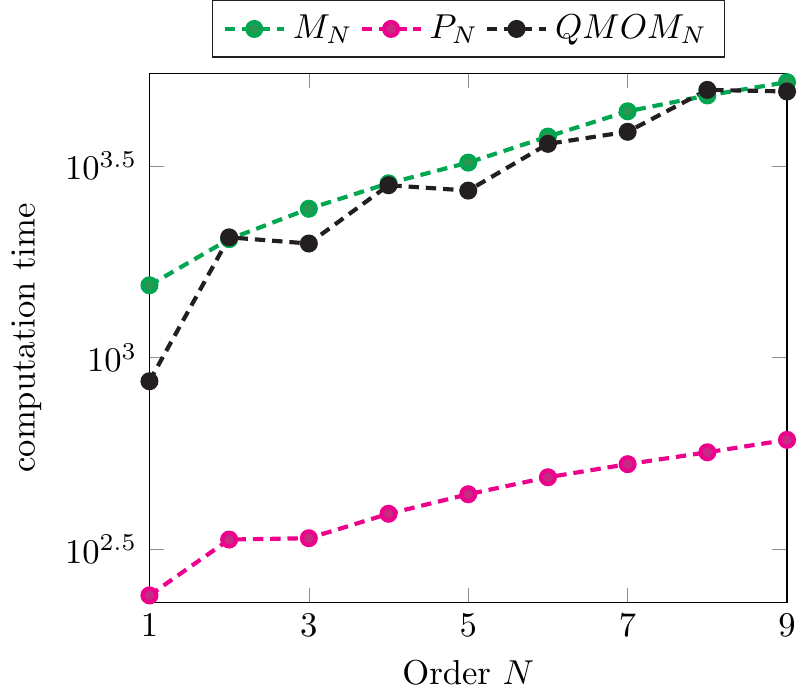}}
		\subfloat[][]{\includegraphics[width=0.45\textwidth]{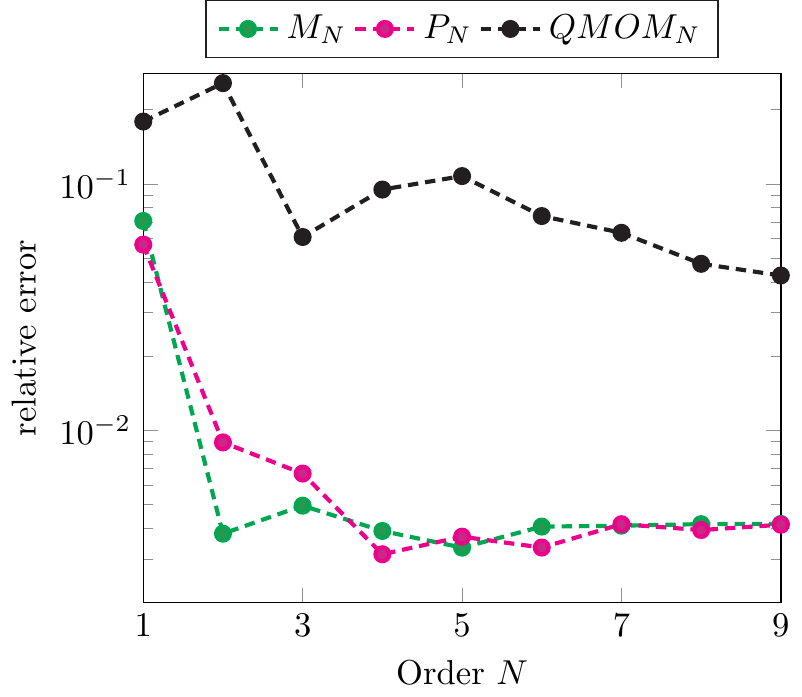}}
	
	\subfloat[][]{\includegraphics[width=0.45\textwidth]{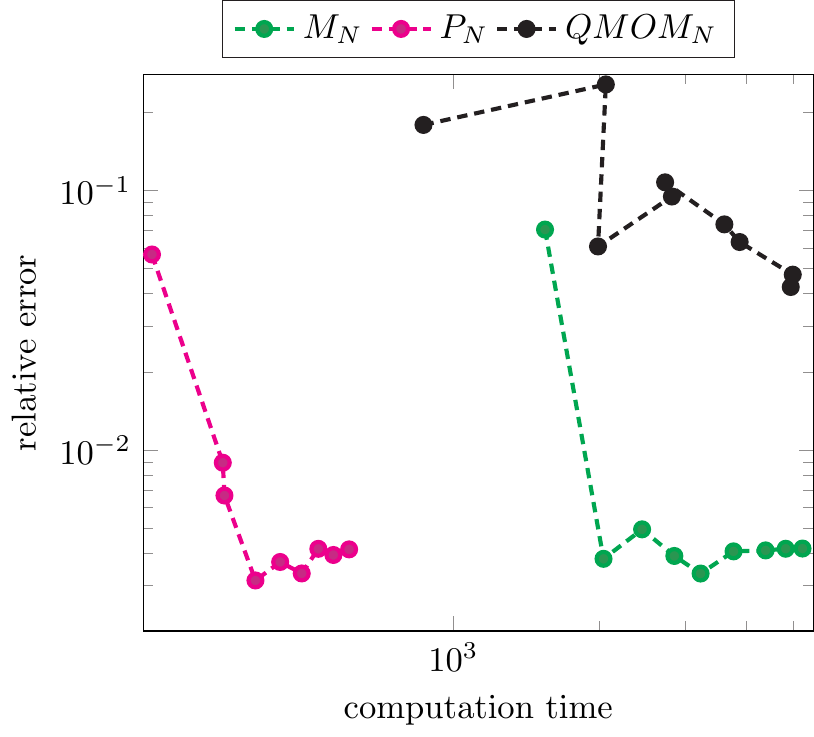}}	
	\caption{Lid driven cavity case: In subfigure (a) the order of moments $N$ is plotted against the computation time, for the maximum entropy model, the polynomial closure and the $QMOM$ and in subfigure (b) the order $N$ is plotted against the relative $E_2$ error for all three methods, in subfigure (c) we can see the computation time plotted against the relative $E_2$ error, again for all three models. In all three subfigures the computations are only done up to the ninth order. The reference solution for the error evaluation is computed by the finite volume scheme (FVS) with $n_v=2000.$}
			\label{figure 2d}
					\end{center}
		\end{figure}

\setcounter{subfigure}{0} 
		\noindent
The comparison of $P_N$ (magenta line), $M_N$ (green line) and $QMOM_N$ (black line) is illustrated in figure \ref{figure 2d}. Again, the densities for all three models are computed with changing values for the order. The numerical solution is computed up to ninth order. The figure is subdivided in the same way as figure \ref{figure break}. It is demonstrated that with increasing order the computation time also increases and the relative $L_2$ error decreases.
Again the Newton algorithm is responsible for the big difference in computation time  between $M_N$ and $P_N.$ 
		\noindent

\begin{rem}
The $P_N$ method is equivalent to the well-known \emph{discrete ordinates} in case of radiative transfer equations \cite{CaseZweifel}, implying that meaningful initial conditions (i.e. those where the underlying $P_N$ distribution is positive for all $v$) will stay meaningful. While this is the case in the examples presented, this cannot be guaranteed in all cases.
This requires to carefully choose the situations where the $P_N$ closure can be applied, if positivity of the kinetic solution is required. If the initial conditions are such that the underlying $P_N$ distribution is not positive, one has to choose $M_N$ or the $QMOM$ model which do not suffer from the above mentioned drawback \cite{Jun00}.
\end{rem}

%% file: Sections/outlook.tex
\section{Conclusions}
\label{outlook}
We have compared a linear moment closure, a nonlinear maximum entropy moment closure and the quadrature method of moments numerically in terms of accuracy, order of the moments and computation times for a population balance equation. 
In the cases where breakage is involved the computation time of the quadrature method of moments and the maximum entropy closure are alike, in contrast to the pure aggregation example. In all cases, the linear closure is the fastest of the three methods. On the other hand, for the linear moment closure  the underlying  distribution function  is not necessarily positive, which can lead to meaningless solutions in some applications. 
The error produced by the linear closure and the maximum entropy closure are smaller than the ones produced by the quadrature method of moments.

Future research should include other versions of $QMOM$, like $EQMOM$ or $SQMOM$ \cite{Attarakih2009,Yuan2012}, investigating if the classical moment models maintain their superiority. Furthermore, low-order ($N\leq 2$) partial moment methods on a partition of the velocity space should be taken into account, simplifying the inversion of the moment problem. These models should be compared with similar low-order $QMOM$ models regarding error and run-time. Additionally, multi-dimensional PBEs and their moment approximations should be investigated. Finally, the comparison of our simulations with real experiments would be interesting.

\section*{Acknowledgements}
Funding by the Deutsche Forschungsgemeinschaft (DFG) within the RTG GrK 1932 "Stochastic Models for Innovations in the Engineering Sciences" and by BMBF Verbundprojekt 05M2016-proMT is gratefully acknowledged.